\documentclass{article}
\usepackage[utf8]{inputenc}

\usepackage{subfigure}
\usepackage{tikz}

\usepackage[utf8]{inputenc}

\usepackage[utf8]{inputenc}
\usepackage[english]{babel}
\usepackage{amsmath}
\usepackage{times}
\usepackage{graphicx}
\usepackage{float}
\usepackage[margin=1in]{geometry}
\usepackage{blindtext}
\usepackage{amssymb}
\usepackage{tikz}
\usepackage{amsthm}
\usepackage{mathtools}
\usepackage{textcmds}

\usepackage{comment}
\usepackage{hyperref}

\newtheorem{thm}{Theorem}[section]
\newtheorem{cor}[thm]{Corollary}
\newtheorem{lem}[thm]{Lemma}

\title{\bf Spectral extremal problems for degenerate graphs\thanks {Research was partially supported by the National
		Nature Science Foundation of China (grant number
		12331012)}}
\date{}
\author {Jiadong Wu$^{1}$, \,  Liying Kang$^{1,2}$\thanks{\em Corresponding author. Email address: lykang@shu.edu.cn (L. Kang), 1753381890@qq.com (J. Wu), 1051466287@qq.com (Z. Ni)}, \, Zhenyu Ni$^{3}$\\
	{\small $^{1}$ Department of Mathematics, Shanghai University,
		Shanghai 200444, P.R. China}\\
	{\small$^{2}$Newtouch Center for Mathematics of Shanghai University,
		Shanghai,  China, 200444}\\
	{\small $^{3}$ Department of Mathematics, Hainan University,
		Haikou 570228, P.R. China}}

\begin{document}
	\maketitle
	
	\begin{abstract}
		A  family of graphs is called degenerate if it contains at least one bipartite graph. In this paper,  we investigate the spectral extremal problems for a degenerate family of graphs $\mathcal{F}$. By employing  covering and independent covering of graphs, we establish a spectral stability result for
		$\mathcal{F}$. Using this stability result, we prove two general theorems that characterize spectral extremal graphs for a broad class of graph families $\mathcal{F}$ and imply several new and known results. Meanwhile, we establish the correlation between extremal graphs and spectral extremal graphs for $\mathcal{F}$.
\bigskip
		
\medskip

\noindent{\bf AMS (2000) subject classification:}  05C50; 05C35		
	\end{abstract}
	
	\section{Introduction}

	For a family of graphs $\mathcal{F}$, a graph is $\mathcal{F}$-free if it does not contain any member of $\mathcal{F}$ as a subgraph. The Tur\'an type problem is to determine the maximum number of edges in an $\mathcal{F}$-free graph of order $n$. The extremal value is called the \textit{Tur\'an number}  of $\mathcal{F}$ and denoted by $\operatorname{ex}(n,\mathcal{F})$. Let $\mathrm{Ex}(n,\mathcal{F})$ be the set of $\mathcal{F}$-free graphs of order $n$ with $\operatorname{ex}(n,\mathcal{F})$ edges.
	Let $\chi(\mathcal{F})=\min_{F\in \mathcal{F}} \chi(F)$. The famous Erd\H{o}s-Stone-Simonovits theorem \cite{ES} states that
	$$\operatorname{ex}(n,\mathcal{F})=\bigg(1-\frac{1}{\chi(\mathcal{F})-1}\bigg)\frac{n^2}{2}+o(n^2).$$
	
	A  family of graphs $\mathcal{F}$ is called \textit{degenerate} if it contains at least one bipartite graph. If $\mathcal{F}$ is a degenerate graph family, then the above problem is called a \textit{ degenerate Tur\'an type problem}. For a degenerate graph family $\mathcal{F}$, the Erd\H{o}s-Stone-Simonovits theorem only gives the bound $\operatorname{ex}(n,\mathcal{F})=o(n^2)$. Over the past decades, there has been a wealth of research on the degenerate Tur\'an type problems. Recently, Alon and Frankl \cite{AF} determined the Tur\'an number of $\{K_{k}, M_{s+1}\}$. Later, Gerbner \cite{G}, Zhu and Chen \cite{ZC} further studied the  value of $\operatorname{ex}(n,\{F, M_{s+1}\})$ for $F$ being a non-bipartite graph or a bipartite graph, respectively. Additionally, Katona and Xiao \cite{KX} proposed an intriguing conjecture on the Tur\'an number of $\{P_k, F\}$, where $F$ is a connected non-bipartite graph. This conjecture was verified by Liu and Kang \cite{LK}.
	For more related results, we refer readers to \cite{FS, FZC, N, XK,Y}.

	
	The spectral radius of $G$, denoted by $\lambda(G)$, is the largest
	eigenvalue of its adjacency matrix $A(G)$. In this paper we investigate spectral analogues of the Tur\'an type problem for graphs, which was proposed by Nikiforov \cite{N1}. Let $\mathcal{F}$ be a family of graphs. The \textit{spectral Tur\'an type problem} is to determine the maximum spectral radius among all $\mathcal{F}$-free graphs of order $n$. Let $\operatorname{spex}(n,\mathcal{F})$ be the maximum value of spectral radius over all $\mathcal{F}$-free graphs of order $n$. The graph attaining $\operatorname{spex}(n,\mathcal{F})$ is called a \textit{spectral extremal graph}. The family of all spectral extremal graphs is denoted by $\mathrm{Ex_{sp}}(n,\mathcal{F})$. There are many results on spectral Tur\'an type problems. In particular, many papers have  investigated the problem of determining when
	$\mathrm{Ex_{sp}}(n,\mathcal{F})\subseteq \mathrm{Ex}(n,\mathcal{F})$ can be guaranteed.
	Wang, Kang and Xue \cite{WKX} confirmed a conjecture of Cioab\u{a}, Desai and Tait \cite{CDT0} and showed that $\mathrm{Ex_{sp}}(n,\mathcal{F})\subseteq \mathrm{Ex}(n,\mathcal{F})$ if $\operatorname{ex}(n,\mathcal{F})=e(T_{r}(n))+O(1)$ and $n$ is sufficiently large. Later, Fang, Tait, and Zhai \cite{FTZ} improved this result for any  graph family $\mathcal{F}$ with $\operatorname{ex}(n,\mathcal{F})<e(T_{r}(n))+\lfloor \frac{n}{2r} \rfloor$. Recently, Byrne \cite{B} further extended this result.
	
	In this paper, we focus on  spectral Tur\'an type problems for degenerate graph families. Over the past decades, numerous scholars have  studied degenerate spectral Tur\'an  problems. Wang, Hou, and Ma \cite{WHM} determined the exact value of $\operatorname{spex}(n,\{M_{s+1}, K_{r+1}\})$ for sufficiently large $n$. Jiang, Yuan, and Zhai \cite{JYZ} characterized graphs in $\mathrm{Ex_{sp}}(n,\{M_{s+1}, F\})$ when $n$ is sufficiently large,  $F$ is a non-bipartite graph. For sufficiently large $n$, Cioab\u{a}, Desai and Tait \cite{CDT} proved that every $n$-vertex graph $G$ with $\lambda(G)\geq \lambda(K_{k}\vee I_{n-k})$ contains all trees of order $2k+2$ unless $G= K_{k}\vee I_{n-k}$. Moreover, if $\lambda(G)\geq \lambda(K_{k}\vee(K_{2}\cup I_{n-k-2}))$, then $G$ contains all trees of order $2k+3$ unless $G= K_{k}\vee(K_{2}\cup I_{n-k-2})$.  Later, Cioab\u{a}, Desai and Tait  \cite{CDT1} showed that $\mathrm{Ex_{sp}}(n,C_{2k+2})=\{K_{k}\vee(K_{2}\cup I_{n-k-2})\}$  for sufficiently large $n$ and $k\geq 2$. Wang, Feng and Lu \cite{WFL}  proved that $\mathrm{Ex_{sp}}(n,\{C_{2k+2}, K_{r+1}\})=\{T_{r-1}(k)\vee I_{n-k}\}$ for sufficiently large $n$ and $k\geq r\geq 2$. For $t\geq 1$ and $k\geq 4$, Zhai, Yuan, and You \cite{ZYY} characterized graphs in $\mathrm{Ex_{sp}}(n, tK_{1,k-1}\cup P_k)$ when $n$ is sufficiently large. Recently, Byrne, Desai and Tait \cite{BDT} established a more general result for a degenerate graph family $\mathcal{F}$ when $\mathrm{Ex}(n,\mathcal{F})$ contains  specific subgraphs and $n$ is sufficiently large.

	\subsection{Notation}
	In this paper, we consider simple and undirected graphs. We first introduce some  notations. Let $G=(V(G),E(G))$ denote a simple graph with vertex set $V(G)$ and edge set $E(G)$.  As customary,  $|G|$ and $e(G)$ represent the number of vertices and edges of a graph $G$, respectively.
	For a vertex $u\in V(G)$, let $N_{G}(u)$ denote the set of neighbors of $u$ in $G$, and $d_{G}(u)=|N_{G}(u)|$ denote its degree in $G$. Further,
	for $U\subseteq V(G)$ and $u\in V(G)$, define $d_{U}(u)=|N_{U}(u)|=|N_{G}(u)\cap U|$. Let $\Delta(G)$ denote the maximum degree of $G$. The $r$-partite Tur\'an graph of order $n$ is denoted by $T_{r}(n)$ and  $I_k$ represents an independent set of $k$ vertices.
	A matching of size of size $s+1$ is denoted by $M_{s+1}$, and the matching number of $G$ is denoted by $\nu (G)$.
	A path with $k$ vertices is denoted by $P_k$, while a cycle with $k$ vertices is denoted by $C_k$. The set of cycles of length at least $k$ is denoted by $C_{\geq k}$.
	For disjoint vertex subsets $U, W\subseteq V(G)$, let $e(U,W)$ denote the number of edges of $G$ between $U$ and $W$.
	The subgraph of $G$ induced by $U$ is denoted by $G[U]$.
	Given two graphs $H$ and $F$, $H\cup F$ denotes the disjoint union of $H$ and $F$, and $H\vee F$ denotes the graph obtained from $H\cup F$ by adding all possible edges between $V(H)$ and $V(F)$.
	
	A \textit{covering} of a graph $G$ is a set of vertices that meets all edges of the graph $G$. An \textit{independent covering} of a graph $G$ is an independent set that meets all edges of $G$.
	Clearly, a graph contains an independent covering if and only if it is a bipartite graph.
	Let $\beta(G)$ denote the minimum size of a covering of $G$, and $\beta'(G)$ the minimum size of an independent covering of $G$.
	If $G$ is not a bipartite graph, we define $\beta'(G)=+\infty$.

	Let $\mathcal{F}$ be a degenerate graph family.    The \textit{covering number} $\beta (\mathcal{F})$ of $\mathcal{F}$ is defined by
	$$\beta(\mathcal{F})=\min \{\beta(F)|F\in \mathcal{F}\}.$$
	The \textit{independent covering number} $\beta' (\mathcal{F})$ of $\mathcal{F}$ is defined by
	$$\beta'(\mathcal{F})=\min \{\beta'(F)|F\in \mathcal{F}\}.$$
	If $\beta'(\mathcal{F})=1$, then $\mathcal{F}$ contains a star $K_{1,t}$ for some $t>0$, which implies that $\operatorname{spex}(n,\mathcal{F})\leq t$.
	For the remainder of this work, we restrict our attention to graph families $\mathcal{F}$ satisfying $\beta'(\mathcal{F})\geq 2$. We introduce a family of graphs. Let $\mathcal{M}(\mathcal{F})=\{F[S]|\ F\in \mathcal{F} \textit{ and S is a covering of F with } |S|<\beta'(\mathcal{F})\}$. Define
	\begin{eqnarray*}
		\mathcal{H}{(\mathcal{F})}=\left\{
		\begin{array}{ll}
			K_{\beta'(\mathcal{F})}, &\mbox{if $\beta'(\mathcal{F})= \beta(\mathcal{F})$},\\[2mm]
			\mathcal{M}(\mathcal{F}), &\mbox{if $\beta'(\mathcal{F})> \beta(\mathcal{F})$}.
		\end{array}
		\right.
	\end{eqnarray*}

	Let $F$ be a graph. We say $F$ is a finite graph if $|F|$ is a fixed constant.
	A  family of graphs  $\mathcal{F}$ is called \textit{finite degenerate} if $\mathcal{F}$ is  \textit{degenerate} and for every $F\in \mathcal{F}$, $|F|$ is bounded by a fixed constant.	A graph family $\mathcal{F}$ is said to be \textit{weak finite degenerate} if there exists a bipartite graph $F\in \mathcal{F}$ with $\beta'(F)=\beta'(\mathcal{F})$ and $|F|$ is bounded by a fixed constant.

	\subsection{Main results}

	Denote $\operatorname{ex}_{G}(n, \mathcal{F})$ as the maximum number of edges in an $\mathcal{F}$-free graph on $n$ vertices with $G$ as a subgraph. We write $\mathrm{Ex}_{G}(n, \mathcal{F})$ for the corresponding extremal graphs. Let
	\begin{eqnarray*}
		\mathcal{G}_{0}(\mathcal{F})&=&\{ T\vee I_{n+1-\beta'(\mathcal{F})}|\ T\in\mathrm{Ex}( \beta'(\mathcal{F})-1 , \mathcal{H}(\mathcal{F}))\},\\[2mm]
		\mathcal{G}(\mathcal{F})&=&\{\mathrm{Ex}_{H}(n,\mathcal{F})|\ H\in \mathcal{G}_{0}(\mathcal{F}) \}.
	\end{eqnarray*}
	We have the following results for a finite degenerate graph family.
	
	\begin{thm}\label{T2}
		Suppose $\mathcal{F}$ is a finite degenerate  family of graphs    with $\beta'(\mathcal{F})\geq 2$ and $\operatorname{ex}(n, \mathcal{F})= O(n)$. Let $H=K_{\beta'(\mathcal{F})-1,n+1-\beta'(\mathcal{F})}$.
		For sufficiently large $n$, if
		$$\operatorname{ex}_{H}(n,\mathcal{F})< e(H)+\left\lfloor \frac{n+1-\beta'(\mathcal{F})}{2}\right\rfloor,$$
		then $\mathrm{Ex_{sp}}(n,\mathcal{F})\subseteq \mathcal{G}(\mathcal{F}).$
	\end{thm}

	\begin{thm}\label{T4}
		Suppose $\mathcal{F}$ is a finite degenerate  family of graphs with $\beta'(\mathcal{F})\geq 2$ and $\operatorname{ex}(n, \mathcal{F})= O(n)$.  Let
		$H=K_{\beta'(\mathcal{F})-1,n+1-\beta'(\mathcal{F})}$ and $H_{1}=T\vee I_{n+1-\beta'(\mathcal{F})}$, where $T\in \mathrm{Ex}(\beta'(\mathcal{F})-1, \mathcal{H}(\mathcal{F}))$.
		For sufficiently large $n$ and some $r\in [\frac{1}{2},\frac{3}{4})$,
		if	$$e(H)+ \left\lfloor \frac{n+1-\beta'(\mathcal{F})}{2}\right\rfloor\leq \operatorname{ex}_{H}(n,\mathcal{F})\leq  e(H)+ r n+O(1)$$
		and $$ \operatorname{ex}_{H}(n,\mathcal{F})-\operatorname{ex}_{H_1}(n,\mathcal{F}) \leq \frac{1}{8(\beta'(\mathcal{F})-1)}n ,$$ then $\mathrm{Ex_{sp}}(n,\mathcal{F})\subseteq \mathcal{G}(\mathcal{F}).$
		
	\end{thm}

	\begin{thm}\label{T21}
		Let $\mathcal{F}$ be a finite degenerate  family of graphs with $\beta'(\mathcal{F})\geq 2$. Suppose $\widetilde{H}$ is a graph on $n+1-\beta'(\mathcal{F})$ vertices with $e(\widetilde{H})\leq rn+O(1)$, where $r\in [0,\frac{3}{4})$. For sufficiently large $n$, if
		$$\{T\vee \widetilde{H}| \ T\in \mathrm{Ex}(\beta'(\mathcal{F})-1, \mathcal{H}(\mathcal{F}))\}\subseteq \mathrm{Ex}(n,\mathcal{F}), $$
		then $\mathrm{Ex_{sp}}(n,\mathcal{F})\subseteq\mathcal{G}(\mathcal{F})\subseteq\mathrm{Ex}(n,\mathcal{F})$.
	\end{thm}

	The graph families investigated above are finite. Next, we extend our study to infinite graph families. Recently, Dou, Hu and Peng \cite{DHP} investigated the Tur\'an number of $\{C_{\geq k}, F\}$. We provide a spectral result for $\{C_{\geq k}, F\}$, where $F$ is an arbitrary finite graph.

	\begin{thm}\label{T3}
		Suppose $k\geq 3$ is an integer and $F$ is an arbitrary finite graph. Let $\mathcal{F}=\{C_{\geq k}, F\}$ and $G$ be an arbitrary graph in $\mathrm{Ex_{sp}}(n,\mathcal{F})$.   If $\beta'(\mathcal{F})=\lfloor \frac{k+1}{2} \rfloor$ and $n$ is sufficiently large, then
		\begin{itemize}
			\item [(1)]
			$G\in  \{T\vee I_{n- \frac{k-1}{2}}| \ T\in\mathrm{Ex}( \frac{k-1}{2}, \mathcal{H}(\mathcal{F})) \} $ when $k$ is odd;
			\item [(2)]
			$G\in \{T\vee I_{n- \frac{k}{2}+1}|\ T\in\mathrm{Ex}( \frac{k}{2}-1, \mathcal{H}(\mathcal{F}))\} $ or  $G\in  \{T\vee (K_2 \cup  I_{n- \frac{k}{2}-1})|\ T\in\mathrm{Ex}( \frac{k}{2}-1, \mathcal{H}(\mathcal{F}))\}$
			when $k$ is even.
		\end{itemize}
	\end{thm}
	
	The remainder of this paper is structured as follows.	In Section 2, we present several applications of our principal findings.
	Section 3 introduces the preliminary lemmas. Section 4 is dedicated to the proofs of Theorems \ref{T2}, \ref{T4}, and \ref{T21}.
	The proof of Theorem \ref{T3} is provided in Section 5.

	\section{Applications of main results}
	
	In this section we give some applications of our main results.

Byrne, Desai and Tait \cite{BDT} proved a general theorem for a family of graphs $\mathcal{F}$ under certain conditions.
We derive the following result from Theorem \ref{T2}, which is a version of Byrne, Desai and Tait's result
when $\mathcal{F}$ is a finite family of graphs.


	\begin{thm}\label{BDT1}
		Let $\mathcal{F}$ be a finite family of graphs. Suppose that $\operatorname{ex}(n,\mathcal{F})=O(n)$, $K_{k+1,\infty}$ is not $\mathcal{F}$-free, and  $\mathrm{Ex}_{K_{k,n-k}}(n,\mathcal{F})\ni G$ for large enough  $n$, where one of the following holds:
		\begin{itemize}
			\item[(a)] $G=K_{k,n-k}$
			
			\item[(b)] $G=K_{k}\vee I_{n-k}$
			
			\item[(c)] $G=K_{k}\vee (K_{2}\cup I_{n-k-2}).$

		\end{itemize}
		If (a), (b), or (c) holds,  then for  large enough $n$, $\mathrm{Ex_{sp}}(n,\mathcal{F})=G$.
		
	\end{thm}
	
	\noindent\textbf{Proof.}
	The assumption that  $K_{k+1,\infty}$ is not $\mathcal{F}$-free and $\mathrm{Ex}_{K_{k,n-k}}(n,\mathcal{F})\ni G$ implies that $\beta'(\mathcal{F})= k+1$.
	If $G=K_{k, n-k}\in \mathrm{Ex}_{K_{k,n-k}}(n,\mathcal{F})$, then $\mathcal{G}(\mathcal{F})= K_{k, n-k}$. By Theorem \ref{T2}, $\mathrm{Ex_{sp}}(n,\mathcal{F})=G$.

	If $G=K_{k}\vee I_{n-k}\in \mathrm{Ex}_{K_{k,n-k}}(n,\mathcal{F})$,
	we claim $\beta'(\mathcal{F})=\beta(\mathcal{F})$. Otherwise, if $\beta(\mathcal{F})<\beta'(\mathcal{F})$, then there exists a graph $F\in \mathcal{F}$ with $\beta(F)=\beta(\mathcal{F})\leq \beta'(\mathcal{F})-1=k$. Thus, $F\subseteq K_{k}\vee I_{n-k}$, which contradicts the fact $K_{k}\vee I_{n-k}$ is $\mathcal{F}$-free.
	Then $\beta'(\mathcal{F})=\beta(\mathcal{F})=k+1$ and it follows that $\mathrm{Ex}(k, \mathcal{H}(\mathcal{F}))=K_{k}$. This together with the fact $ K_{k}\vee I_{n-k}\in \mathrm{Ex}_{K_{k,n-k}}(n,\mathcal{F})$ implies that $\mathcal{G}(\mathcal{F})= \{K_{k}\vee I_{n-k}\}$.  According to Theorem \ref{T2}, $\mathrm{Ex_{sp}}(n,\mathcal{F})=G$.
	Similarly, if $G=K_{k}\vee (K_{2}\cup I_{n-k-2})\in \mathrm{Ex}_{K_{k,n-k}}(n,\mathcal{F})$,  we can deduce that $\mathrm{Ex_{sp}}(n,\mathcal{F})=G$. \qed
	
\noindent	{\bf Remark:}  As special cases of Theorem \ref{BDT1}, we recovered results from \cite{CLZ2019, PC2011, CDT, FLSZ2024, FYZ2007, N2010}.

	Using Theorem \ref{T2},  Theorem \ref{T4}, we get the following result due to Byrne, Desai and Tait \cite{BDT}.
	
	\begin{thm}[Byrne, Desai and Tait \cite{BDT}]\label{BDT}
		Let $\mathcal{F}$ be a family of graphs. Suppose that $\operatorname{ex}(n,\mathcal{F})=O(n)$, $K_{k+1,\infty}$ is not $\mathcal{F}$-free, and for $n$ large enough $\mathrm{Ex}_{K_{k,n-k}}(n,\mathcal{F})\ni K_{k}\vee X$, where $e(X)\leq rn+O(1)$ for some $r\in[0,3/4)$ and $\mathcal{F}$ is finite.
		Then for $n$ large enough, $\mathrm{Ex_{sp}}(n,\mathcal{F})\subseteq\mathrm{Ex}_{K_{k}\vee I_{n-k}}(n,\mathcal{F}).$
	\end{thm}
	
	\noindent\textbf{Proof.}
	Since $K_{k}\vee X\in \mathrm{Ex}_{K_{k,n-k}}(n,\mathcal{F})$ and $K_{k+1,\infty}$ is not $\mathcal{F}$-free, similar as the proof of Theorem \ref{BDT1}, one can easily verify that $\beta'(\mathcal{F})=\beta(\mathcal{F})= k+1$.
	
	Let $H=K_{k,n-k}, H_1=K_{k}\vee I_{n-k}$. The assumption $K_{k}\vee X\in \mathrm{Ex}_{K_{k,n-k}}(n,\mathcal{F})$ implies that $\operatorname{ex}_{H_1}(n,\mathcal{F})= \operatorname{ex}_{H}(n,\mathcal{F})$ and
	$\operatorname{ex}_{H}(n,\mathcal{F})= e(H)+ e(X)+\binom{k}{2}$.
	If $e(X)<\left\lfloor \frac{n-k}{2}\right\rfloor-\binom{k}{2}$, then $\operatorname{ex}_{H}(n,\mathcal{F})< e(H)+\lfloor \frac{n-k}{2}\rfloor.$  By Theorem \ref{T2}, $\mathrm{Ex_{sp}}(n,\mathcal{F})\subseteq\mathrm{Ex}_{K_{k}\vee I_{n-k}}(n,\mathcal{F})$. If $\left\lfloor \frac{n-k}{2}\right\rfloor-\binom{k}{2}\leq e(X)\leq rn+O(1)$ for some $r\in [\frac{1}{2},\frac{3}{4})$, then
	$$e(H)+\left \lfloor \frac{n-k}{2} \right\rfloor\leq \operatorname{ex}_{H_1}(n,\mathcal{F})= \operatorname{ex}_{H}(n,\mathcal{F})\leq  e(H)+rn+O(1).$$
	By Theorem \ref{T4}, $\mathrm{Ex_{sp}}(n,\mathcal{F})\subseteq\mathrm{Ex}_{K_{k}\vee I_{n-k}}(n,\mathcal{F}).$
	\qed

	Next we consider the finite degenerate  family of graphs which contains a matching. Suppose  $\mathcal{F}=\{M_{s+1},F\}$, where $F$ is an arbitrary finite graph. For sufficiently large $n$, Zhu and Chen \cite{ZC} proved that $\mathrm{Ex}(n, \mathcal{F})=\{T\vee I_{n-s}|\ T\in\mathrm{Ex}(s, \mathcal{H}(\mathcal{F}))\}$ if $\beta'(\mathcal{F})=s+1$ and $\operatorname{ex}(n, \mathcal{F})=(\beta'(F)-1)n+O(1)$ if $\beta'(\mathcal{F})\leq s$. We give a spectral version of their results.
	
	\begin{thm}\label{A1}
		Let  $\mathcal{F}=\{M_{s+1},F\}$, where  $F$ is an arbitrary finite graph. For sufficiently large $n$, the following statements hold.
		\begin{itemize}
			\item [(1)] If $\beta'(\mathcal{F})=s+1$, then $\mathrm{Ex_{sp}}(n,\mathcal{F})\subseteq \mathrm{Ex}(n, \mathcal{F}).$
			\item [(2)] If $2\leq \beta'(\mathcal{F})\leq s$, then $\mathrm{Ex_{sp}}(n,\mathcal{F})\subseteq \mathcal{G}(\mathcal{F}).$
		\end{itemize}
	\end{thm}
	\noindent\textbf{Proof.} (1)
	If $\beta'(\mathcal{F})=s+1$,  by combining Theorem  \ref{T21} and Zhu and Chen's  result \cite{ZC} that $\mathrm{Ex}(n, \mathcal{F})=\{T\vee I_{n-s}|\ T\in\mathrm{Ex}(s, \mathcal{H}(\mathcal{F}))\}$,
	we conclude that $\mathrm{Ex_{sp}}(n,\mathcal{F})\subseteq \mathrm{Ex}(n, \mathcal{F})$.

	(2)	If $2\leq \beta'(\mathcal{F})\leq s$, then $\beta'(\mathcal{F})=\beta'(F)$. Let $H=K_{\beta'(F)-1, n+1-\beta'(F)}$.
	Assume $G$ is an  $\mathcal{F}$-free graph on $n$ vertices with $H\subseteq G$.
	Let $A$ and $B$ be the  color classes
	of $H$ with $|A|=\beta'(F)-1$ and $|B|=n+1-\beta'(F)$.
	For any vertex $u\in B$, it follows that $d_{B}(u)\leq |F|-1$ since $F\subseteq K_{\beta'(F), |F|}$. On the other hand, $\nu (G[B])\leq s$ since $G$ is $M_{s+1}$-free. By the  result of Chv\'atal and Hanson \cite{CH},
	$e(G[B])\leq \nu (G[B])(\Delta(G[B])+1)\leq s|F|$. Thus,  $e(G)\leq e(H)+s|F|+\binom{\beta'(F)-1}{2}$. Let $c=s|F|+\binom{\beta'(F)-1}{2}$. So $\operatorname{ex}_{H}(n,\mathcal{F})\leq e(H)+c$. By Theorem \ref{T2}, $\mathrm{Ex_{sp}}(n,\mathcal{F})\subseteq \mathcal{G}(\mathcal{F}).$ This completes the proof. \qed
	
	Denote $G(n,r,s)$ as the complete $r$-partite graph on $n$ vertices with one part of order $n-s$ and each other part of order $\left \lfloor \frac{s}{r-1}\right \rfloor$ or $\left \lceil \frac{s}{r-1} \right \rceil$. Recently, Alon and Frankl \cite{AF} established that $\mathrm{Ex}(n,\{M_{s+1}, K_{r+1}\})=G(n,r,s)$ for sufficiently large $n$. Clearly, for $r\geq 2$, $\beta'(\{M_{s+1}, K_{r+1}\})=s+1$. Then, by Theorem \ref{A1},  the following result of Wang, Hou and Ma \cite{WHM} follows immediately.

	\begin{cor}[Wang, Hou and Ma \cite{WHM}]
		Let $r\geq 2$ and $s\geq 1$ be two integers. For sufficiently large $n$,
		$\mathrm{Ex_{sp}}(n,\{M_{s+1}, K_{r+1}\})=G(n,r,s)$.
	\end{cor}
	
	A graph is called \textit{edge-critical} if it contains an edge whose deletion decreases its chromatic number. Let $F$ be an edge-critical graph with $\chi(F)=r+1\geq3$ and $\mathcal{F}=\{M_{s+1}, F\}$. For large enough $n$, Alon and Frankl \cite{AF} proved that $\operatorname{Ex}(n,\mathcal{F})=G(n,r,s)$. Clearly, $\beta'(\mathcal{F})=s+1$.  By applying Theorem \ref{A1}, we derive a spectral result for  $\mathcal{F}$.

	\begin{cor}
		Let $F$ be an edge-critical graph with $\chi(F)=r+1\geq 3$ and $\mathcal{F}=\{M_{s+1}, F\}$. If $n$ is  sufficiently large, then
		$\mathrm{Ex_{sp}}(n,\mathcal{F})=G(n,r,s).$
	\end{cor}
	
	Let $\mathcal{J}$ be a finite degenerate  family of graphs with $\beta'(\mathcal{J})=\beta(\mathcal{J})=s\geq 2$ and $\operatorname{ex}(n,\mathcal{J})=O(n)$. Suppose there exists a constant $c\geq 1$ such that $\operatorname{ex}_{H}(n,\mathcal{J})\leq e(H)+c$, where $H=K_{s-1,n+1-s}$. Let $r\geq 2$ and $\mathcal{F}=\{K_{r+1},\mathcal{J}\}$.
	Obviously, $\beta(\mathcal{F})=\min \{r,s\}$ and  $\beta'(\mathcal{F})=\beta'(\mathcal{J})=s\geq 2$.
	One can easily check that $\mathcal{G}(\mathcal{F})=\{T_{r-1}(s-1)\vee I_{n-s+1}\}$ if $r\leq s$, and  $\mathcal{G}(\mathcal{F})=\mathrm{Ex}_{T}(n,\mathcal{F})$, where $T=K_{s-1}\vee I_{n-s+1}$, if $r\geq s+1$. Therefore, by Theorem \ref{T2}, we obtain the following result.
	
	\begin{cor}\label{A2}
		Let $r\geq 2$, $\mathcal{J}$ be a finite degenerate graph family with $\beta'(\mathcal{J})=\beta(\mathcal{J})\geq 2$ and $\operatorname{ex}(n,\mathcal{J})=O(n)$. Assume there exists a constant $c\geq 1$ such that $\operatorname{ex}_{H}(n,\mathcal{J})\leq e(H)+c$, where $H=K_{\beta'(\mathcal{J})-1,n+1-\beta'(\mathcal{J})}$. Let $\mathcal{F}=\{K_{r+1},\mathcal{J}\}$.  For sufficiently large $n$, the following statements hold.
		\begin{itemize}
			\item[(1)] If $r\leq \beta'(\mathcal{J})$, then $\mathrm{Ex_{sp}}(n,\mathcal{F})= T_{r-1}(\beta'(\mathcal{J})-1)\vee I_{n-\beta'(\mathcal{J})+1}$.
			\item [(2)] If $r\geq \beta'(\mathcal{J})+1$, then $\mathrm{Ex_{sp}}(n,\mathcal{F})\subseteq \mathrm{Ex}_{T}(n,\mathcal{F})$, where $T=K_{\beta'(\mathcal{J})-1}\vee I_{n-\beta'(\mathcal{J})+1}$.
		\end{itemize}
	\end{cor}
	
	The famous Erd\H{o}s-S\'os conjecture states that every graph  with average degree larger than $t-2$ contains all trees of order $t$. By Corollary \ref{A2}, we can deduce the following result concerning spectral Erd\H{o}s-S\'os theorem due to Wang, Feng and Lu \cite{WFL}.
	
	\begin{thm}[Wang, Feng and Lu \cite{WFL}]\label{B1}
		Let $t$ and $r$ be two integers such that $t\geq r\geq2$. Let $G$ be a $K_{r+1}$-free graph with order $n$. For sufficiently large $n$, if $\lambda(G)\geq \lambda(T_{r-1}(t)\vee I_{n- t})$, then $G$ contains all trees on $2t+2$ or $2t+3$ vertices unless $G=T_{r-1}(t)\vee I_{n-t}$.
	\end{thm}
	
	\noindent\textbf{Proof.}
	Let $\mathcal{F}$ be the set of all finite graphs that contain all trees on $2t+2$ vertices. The famous Cayley formula shows that the number of trees of order $2t+2$ is at most $(2t+2)^{2t}$. Note that $G$ contains all trees of order $2t+2$ if and only if $G$ is not $\mathcal{F}$-free. Next, we show that $\beta'(\mathcal{F})=\beta(\mathcal{F})=t+1$. Since $K_{t+1,2t+2}\in \mathcal{F}$, we have $\beta'(\mathcal{F})\leq t+1$. We claim $\beta(\mathcal{F})\geq t+1$. Otherwise, suppose  there exists a graph $F\in \mathcal{F}$ such that $\beta(F)\leq t$. Then $F\subseteq K_{t}\vee I_{|F|}$. Note that $K_{t}\vee I_{|F|}$ is $P_{2t+2}$-free, which contradicts the fact that $F$ contains all trees on $2t+2$ vertices. Thus, $t+1\leq \beta(\mathcal{F})\leq \beta'(\mathcal{F})\leq t+1$, which proves that $\beta'(\mathcal{F})=\beta(\mathcal{F})=t+1$.

	Let $T_t$ be an arbitrary tree of order $t$. Cioab\u{a}, Desai and Tait \cite{CDT} showed that for any $n$,  $\frac{1}{2}(t-2)n\leq \operatorname{ex}(n,T_t)\leq (t-2)n$.
	In the same paper, they further established  that for $H=K_{t,n-t}$, the inequality    $\operatorname{ex}_{H}(n,\mathcal{F})\leq e(H)+\binom{t}{2}$ holds.
	Since $\operatorname{ex}(n,T_{2t+2})\leq 2tn$, it follows that  $\operatorname{ex}(n,\mathcal{F})=O(n)$.  Then, given $\beta'(\mathcal{F})>t\geq r\geq2$,  part (1) of Corollary \ref{A2} yields
	$\mathrm{Ex_{sp}}(n,\{K_{r+1},\mathcal{F}\})= T_{r-1}(t)\vee I_{n-t}$. Recall that $G$ is $K_{r+1}$-free. Therefore, if $\lambda(G)\geq \lambda(T_{r-1}(t)\vee I_{n- t})$, then $G$ contains all trees of order $2t+2$ unless $G=T_{r-1}(t)\vee I_{n-t}$. The proof for trees with $2t+3$ vertices is analogous, so we omit the details for brevity. The proof is complete. \qed

	Recently, Katona and Xiao \cite{KX} determined $\operatorname{ex}(n,\{K_{r+1},P_{k+1}\})$ for $k\geq 2r+1$ and sufficiently large $n$. We give a spectral version of the result.  Set $\mathcal{J}=\{P_{k+1}\}$ in Corollary \ref{A2}, we obtain the following result.
	
	\begin{cor}\label{A3}
		Let $r\geq 2$ and $k\geq 3$ be two integers.  For   sufficiently large $n$, the following statements hold.
		\begin{enumerate}
			\item[(1)]  If  $ r\leq \lfloor \frac{k+1}{2} \rfloor $, then $\mathrm{Ex_{sp}}(n, \{K_{r+1}, P_{k+1}\}) = T_{r-1}(\lfloor \frac{k-1}{2} \rfloor) \vee I_{n -\lfloor \frac{k-1}{2} \rfloor}$.
			
			\item[(2)] If $ r\geq \lfloor \frac{k+1}{2} \rfloor +1 $ and $k$ is odd, then $\mathrm{Ex_{sp}}(n, \{K_{r+1}, P_{k+1}\}) = K_{\frac{k-1}{2}} \vee I _{n - \frac{k-1}{2}}$. \\
			If $ r\geq \lfloor \frac{k+1}{2} \rfloor +1 $ and $k$ is even, then $\mathrm{Ex_{sp}}(n, \{K_{r+1}, P_{k+1}\}) = K_{\frac{k}{2} - 1} \vee \big( K_2 \cup I_{n - \frac{k}{2} - 1} \big)$.
		\end{enumerate}
	\end{cor}

	A \textit{linear forest} is a graph whose connected components are paths. Denote $\mathcal{L}_k$ as the family of all linear forests with size $k$. Note that $P_{k+1}\in \mathcal{L}_k$. One can easily check that all graphs in $\operatorname{EX_{sp}}(n, \{K_{r+1}, P_{k+1}\})$ are also $\{K_{r+1}, \mathcal{L}_k\}$-free. Thus, $\mathrm{Ex_{sp}}(n, \{K_{r+1}, \mathcal{L}_k\})=\mathrm{Ex_{sp}}(n, \{K_{r+1}, P_{k+1}\})$. Then we can get the result of Zhai and Yuan \cite{ZY} by Corollary \ref{A3}.
	
	\begin{cor} [Zhai and Yuan \cite{ZY}]
		Let $r\geq 2$ and $k\geq 3$ be two integers. For   sufficiently large $n$, the following statements hold.
		\begin{enumerate}
			\item[(1)]  If  $ r\leq \lfloor \frac{k+1}{2} \rfloor $, then $\mathrm{Ex_{sp}}(n, \{K_{r+1}, \mathcal{L}_k\}) = T_{r-1}(\lfloor \frac{k-1}{2} \rfloor) \vee I_{n -\lfloor \frac{k-1}{2} \rfloor}$.
			
			\item[(2)] If $ r\geq \lfloor \frac{k+1}{2} \rfloor +1 $ and $k$ is odd, then $\mathrm{Ex_{sp}}(n, \{K_{r+1}, \mathcal{L}_k\}) = K_{\frac{k-1}{2}} \vee I _{n - \frac{k-1}{2}}$. \\
			If $ r\geq \lfloor \frac{k+1}{2} \rfloor +1 $ and $k$ is even, then $\mathrm{Ex_{sp}}(n, \{K_{r+1}, \mathcal{L}_k\}) = K_{\frac{k}{2} - 1} \vee \big( K_2 \cup I_{n - \frac{k}{2} - 1} \big)$.
		\end{enumerate}
	\end{cor}

	The  families of graphs we have investigated above are all finite. We now turn to studying some infinite  families of graphs.
	Recently, Dou, Ning, and Peng \cite{DNP} determined the  Tur\'an number of $\{C_{\geq k}, K_{r+1}\}$.
	Let $\mathcal{F}=\{C_{\geq k}, K_{r+1}\}$. One can easily verify that $\beta(\mathcal{F})=\mbox{min} \{r, \lfloor \frac{k+1}{2} \rfloor\}$ and
	$\beta'(\mathcal{F})= \lfloor \frac{k+1}{2} \rfloor$. If  $r\leq \lfloor \frac{k+1}{2} \rfloor$, then $\mathcal{H}(\mathcal{F})=\{K_r\}$.
	If $r\geq \lfloor \frac{k+1}{2} \rfloor+1$, then $\beta(\mathcal{F})=\beta'(\mathcal{F})$ and
	$\mathcal{H}(\mathcal{F})=\left\{K_{\lfloor \frac{k+1}{2} \rfloor}\right\}$.
	By Theorem \ref{T3}, we derive a spectral result for $\{C_{\geq k}, K_{r+1}\}$.

	\begin{thm}\label{A4}
		Let $r\geq 2$ and $k\geq 3$ be two integers. For   sufficiently large $n$, the following statements hold.
		\begin{itemize}
			\item[(1)]  If  $ r\leq \lfloor \frac{k+1}{2} \rfloor $, then $\mathrm{Ex_{sp}}(n, \{C_{\geq k}, K_{r+1}\}) = T_{r-1}(\lfloor \frac{k-1}{2} \rfloor) \vee I_{n -\lfloor \frac{k-1}{2} \rfloor}$.
			
			\item[(2)] If $ r\geq \lfloor \frac{k+1}{2} \rfloor +1 $ and $k$ is odd, then $\mathrm{Ex_{sp}}(n, \{C_{\geq k}, K_{r+1}\}) = K_{\frac{k-1}{2}} \vee I _{n - \frac{k-1}{2}}$. \\
			If $ r\geq \lfloor \frac{k+1}{2} \rfloor +1 $ and $k$ is even, then $\mathrm{Ex_{sp}}(n, \{C_{\geq k}, K_{r+1}\}) = K_{\frac{k}{2} - 1} \vee \big( K_2 \cup I_{n - \frac{k}{2} - 1} \big)$.
		\end{itemize}
	\end{thm}

	Recently, Zhao and Lu \cite{ZL} studied the Tur\'an number of $\{C_{\geq k}, M_{s+1}\}$ when $n$ is large enough and $s\geq \lfloor \frac{k-1}{2} \rfloor$. Let $\mathcal{F}=\{C_{\geq k}, M_{s+1}\}$. If $s\geq \lfloor \frac{k-1}{2}\rfloor$, then $\beta(\mathcal{F})=\beta'(\mathcal{F})=\lfloor \frac{k+1}{2} \rfloor$. Thus $\mathcal{H}(\mathcal{F})=\left \{K_{\lfloor \frac{k+1}{2} \rfloor} \right\}$ and $\mathrm{Ex}(\lfloor \frac{k-1}{2} \rfloor,\mathcal{H}(\mathcal{F}))=K_{\lfloor \frac{k-1}{2} \rfloor}$. By Theorem \ref{T3}, we  get a spectral result for $\{C_{\geq k}, M_{s+1}\}$.

	\begin{cor}
		Let $k\geq 3$ and $s\geq \lfloor \frac{k-1}{2}\rfloor$ be two integers. For   sufficiently large $n$, the following statements hold.
		\begin{itemize}
			\item [(1)] If $k$ is odd, then $\mathrm{Ex_{sp}}(n, \{C_{\geq k},M_{s+1}\}) = K_{\frac{k-1}{2}} \vee I _{n - \frac{k-1}{2}}$.
			
			\item [(2)] If $k$ is even and $s= \frac{k}{2}-1$, then $\mathrm{Ex_{sp}}(n, \{C_{\geq k}, M_{s+1}\}) = K_{\frac{k}{2} - 1} \vee I_{n - \frac{k}{2} + 1}$.\\
			if $k$ is even and $s\geq \frac{k}{2}$, then $\mathrm{Ex_{sp}}(n, \{C_{\geq k}, M_{s+1}\}) =K_{\frac{k}{2} - 1} \vee \big( K_2 \cup I_{n - \frac{k}{2} - 1} \big)$.
		\end{itemize}
	\end{cor}
	
	\section{Preliminaries}	
	We introduce the following lemmas that will be employed to prove our results.
	
	\begin{lem}[Wu, Xiao and Hong \cite{WXH}]\label{WXH}
		Let $G$ be a connected graph and $\mathbf{x}$ be the Perron vector of $G$. Assume that $u,v$ are two vertices of $G$ with $x_u \geq x_v$ and $\{v_i \mid 1 \leq i \leq s\} \subseteq N_G(v) \setminus (N_G(u) \cup \{u\})$. Let $G' = G - \{vv_i \mid 1 \leq i \leq s\} + \{uv_i \mid 1 \leq i \leq s\}$. Then, $\lambda(G') > \lambda(G)$.
	\end{lem}
	\begin{lem}\label{L0}
		Let $\mathcal{F}$ be a family of graphs. For any graph $G\in \{T\vee I_{n+1-\beta'(\mathcal{F})}|\ T\in \mathrm{Ex}(\beta'(\mathcal{F})-1, \mathcal{H}(\mathcal{F}))\}$,  $G$ is $\mathcal{F}$-free.
	\end{lem}
	
	\noindent\textbf{Proof.} Assume $G=T\vee I_{n+1-\beta'(\mathcal{F})}$.
	It is sufficient to prove that $G$ is $F$-free for any graph  $F\in \mathcal{F}$.
	If $\beta(F)
	\geq \beta'(\mathcal{F})$, we claim that $K_{\beta'(\mathcal{F})-1} \vee I_{n+1-\beta'(\mathcal{F})}$ is $F$-free. Otherwise,  $F\subseteq K_{\beta'(\mathcal{F})-1} \vee I_{n+1-\beta'(\mathcal{F})}$, then $\beta(F)\leq \beta'(\mathcal{F})-1\leq \beta(F)-1$, a contradiction. Thus, $K_{\beta'(\mathcal{F})-1} \vee I_{n+1-\beta'(\mathcal{F})}$ is $F$-free. This together with the fact $G\subseteq K_{\beta'(\mathcal{F})-1} \vee I_{n+1-\beta'(\mathcal{F})}$ implies that $G$ is $F$-free.
	
	If $\beta(F)<\beta'(\mathcal{F})$, then $\beta(\mathcal{F})<\beta'(\mathcal{F})$, so $\mathcal{H}(\mathcal{F})=\mathcal{M}(\mathcal{F})$. If $F\subseteq G$, then there exists a covering $S$ of $F$ with $|S|\leq \beta'(\mathcal{F})-1$ such that $F[S]$ is a subgraph of $T$. This contradicts the fact that all graphs in $\mathrm{Ex}(\beta'(\mathcal{F})-1, \mathcal{M}(\mathcal{F}))$ are $F[S]$-free. Consequently, $G$ is $F$-free.\qed

	We establish the following spectral stability result for a weakly finite degenerate graph family.
	\begin{thm}\label{T1}
		Let $C\geq1$ be a real number, and $\mathcal{F}$ be a weak finite degenerate graph family with $\beta'(\mathcal{F})\geq 2$ and $\operatorname{ex}(n, \mathcal{F})\leq Cn$. Suppose $G$ is an $\mathcal{F}$-free graph with order $n$, and $\mathbf{x}$ is a non-negative eigenvector of $A(G)$ corresponding to $\lambda(G)$, with maximum entry equal to $1$. For sufficiently large $n$ and a sufficiently small $\varepsilon$ satisfying $0<\varepsilon<\big(\frac{1}{32C\beta'(\mathcal{F})}\big)^{40}$, the following statements hold:
		\begin{enumerate}
			\item[(1)]
			If
			$ \lambda(G)\geq \sqrt{(\beta'(\mathcal{F})-1)(n+1-\beta'(\mathcal{F}))}-\sqrt{\varepsilon n} ,$
			then there exists a  $t>(1-2\beta'(\mathcal{F})\varepsilon^{\frac{1}{10}})n$ such that $K_{\beta'(\mathcal{F})-1,t}\subseteq G$.
			\item[(2)] Let $L'$ be the color class of $K_{\beta'(\mathcal{F})-1,t}$ with size $\beta'(\mathcal{F})-1$. Then $x_v\geq 1-\varepsilon^{\frac{1}{10}}$ for any $v\in L'$ and  $x_v\leq \varepsilon^{\frac{1}{8}}$ for any $v\notin L'$.
		\end{enumerate}
	\end{thm}
	\noindent\textbf{Proof of Theorem \ref{T1}.} (1).
	Let $F_{0}$ be a bipartite graph in $ \mathcal{F}$ with $\beta'(F_0)=\beta'(\mathcal{F})$ and $|F_0|\leq l$. Set $\lambda=\lambda(G)$, $\beta'=\beta'(\mathcal{F})$.
	For any $0<\eta\leq 1$, define $L^{\eta}=\{u\in V(G)|\ x_{u}\geq \eta\}$. For a vertex $u\in V(G)$ and a positive integer $i$, denote $N_{i}(u)$ as the set of vertices at distance $i$ from $u$ in $G$.

	
	\begin{lem}\label{L1}
		There exists a constant $D(\eta, \varepsilon)$ depending on $\eta$ and $\varepsilon$ such that $|L^{\eta}|\leq D(\eta,\varepsilon)\sqrt{n}$.
	\end{lem}
	
	\noindent\textbf{Proof.} For any $u\in L^{\eta}$, we have
	
	$$(\sqrt{(\beta'-1)(n+1-\beta')}-\sqrt{\varepsilon n})\eta\leq \lambda x_{u}=\sum_{v \in N_{G}(u)}x_{v}\leq d_{G}(u).$$
	Summing  over all vertices $u\in  L^{\eta}$, we obtain
	$$|L^{\eta}| (\sqrt{(\beta'-1)(n+1-\beta')}-\sqrt{\varepsilon n})\eta\leq \sum_{u\in L^{\eta}}d_{G}(u)\leq \sum_{u\in V(G)}d_{G}(u)\leq 2Cn,$$
	where the last inequality holds since $G$ is $\mathcal{F}$-free. By algebraic manipulation, it follows that
	$$| L^{\eta}|\leq \frac{2C}{\eta \sqrt{\beta'-1-\varepsilon^{\frac{1}{4}}}}\sqrt{n}=D(\eta,\varepsilon)\sqrt{n}.$$ \qed
	
	Let $u$ be an arbitrary vertex in $V(G)$. Denote $L_{i}^{\eta}(u)=N_{i}(u)\cap L^{\eta}$ and $\overline{L_{i}^{\eta}}(u)=N_{i}(u)\backslash L^{\eta}$.
	
	\begin{lem}\label{M5}
		For any vertex $u\in V(G)$, we have $$ (\beta'-1-\varepsilon^{\frac{1}{4}})nx_{u} \leq  d_{G}(u)x_{u}+\sum_{v \in \overline{L_1^\eta}(u)} \sum_{w \in L_1^\eta(u) \cup L_2^\eta(u)} x_w+3C\eta n.$$
	\end{lem}
	
	\noindent\textbf{Proof.}
	For any vertex $u\in V(G)$, we have
	\begin{eqnarray}\label{eq1}
		(\beta'-1-\varepsilon^{\frac{1}{4}})nx_{u} &\leq &\lambda^{2}x_{u}\nonumber\\
		&=&d_{G}(u)x_{u}+\sum_{v\in N_1(u)}\sum_{w\in N_{1}(v)\backslash\{u\}}x_{w}\nonumber\\
		&\leq & d_G(u)x_u + \sum_{v \in N_1(u)} \sum_{w \in L_1^\eta(u) \cup L_2^\eta(u)} x_w + \sum_{v \in N_1(u)} \sum_{w \in \overline{L_1^\eta}(u) \cup \overline{L_2^\eta}(u)} x_w,
	\end{eqnarray}
	where the last inequality holds since $N_{1}(v)\backslash\{u\}\subseteq N_{1}(u)\cup N_{2}(u)= L_1^\eta(u) \cup L_2^\eta(u) \cup  \overline{L_1^\eta}(u) \cup \overline{L_2^\eta}(u)$.

	Since $N_{1}(u)=L_1^\eta(u) \cup \overline{L_1^\eta}(u)$, it follows that
	\begin{eqnarray}\label{eq2}
		\sum_{v \in N_1(u)} \sum_{w \in \overline{L_1^\eta}(u) \cup \overline{L_2^\eta}(u)} x_w &=& \sum_{v \in  L_1^\eta(u)} \sum_{w \in \overline{L_1^\eta}(u) \cup \overline{L_2^\eta}(u)} x_w + \sum_{v \in \overline{L_1^\eta}(u)} \sum_{w \in \overline{L_1^\eta}(u) \cup \overline{L_2^\eta}(u)} x_w \nonumber\\[2mm]
		&\leq &\bigg(e(L_1^\eta(u),  \overline{L_1^\eta}(u) \cup \overline{L_2^\eta}(u))+2e(\overline{L_1^\eta}(u)) +e(\overline{L_1^\eta}(u),\overline{L_2^\eta}(u))\bigg)\eta\nonumber\\[2mm]
		&\leq &2e(G)\eta\leq 2C\eta n.
	\end{eqnarray}
	Similarly,	
	\begin{eqnarray}\label{eq3}
		\sum_{v \in N_1(u)} \sum_{w \in L_1^\eta(u) \cup L_2^\eta(u)} x_w &= &\sum_{v \in  L_1^\eta(u)} \sum_{w \in L_1^\eta(u) \cup L_2^\eta(u)} x_w + \sum_{v \in \overline{L_1^\eta}(u)} \sum_{w \in L_1^\eta(u) \cup L_2^\eta(u)} x_w \nonumber\\[2mm]
		&\leq & 2e(L_1^\eta(u))+e(L_1^\eta(u),L_2^\eta(u))+\sum_{v \in \overline{L_1^\eta}(u)} \sum_{w \in L_1^\eta(u) \cup L_2^\eta(u)} x_w.
	\end{eqnarray}
	Note that $L_1^\eta(u) \cup L_2^\eta(u)\subseteq L^\eta$ and $G$ is $\mathcal{F}$-free. Combining this with Lemma \ref{L1}, we get
	\begin{eqnarray}\label{eq4}
		2e(L_1^\eta(u))+e(L_1^\eta(u),L_2^\eta(u))\leq 2C|L_1^\eta(u)|+C(|L_1^\eta(u)|+|L_2^\eta(u)|)
		\leq 3C|L^\eta|\leq 3CD(\eta,\varepsilon)\sqrt{n}.
	\end{eqnarray}

	\noindent
	By (\ref{eq1})-(\ref{eq4}) and the assumption that $n$ is large enough, we obtain
	
	\begin{equation*}
		\begin{aligned}
			(\beta'-1-\varepsilon^{\frac{1}{4}})nx_{u}
			\leq d_{G}(u)x_{u}+\sum_{v \in \overline{L_1^\eta}(u)} \sum_{w \in L_1^\eta(u) \cup L_2^\eta(u)} x_w+3C\eta n.
		\end{aligned}
	\end{equation*}
	\qed
	
	\begin{lem}\label{L2}
		Let $m(\varepsilon)=\frac{4C^2}{(\beta'-1-2\varepsilon^{\frac{1}{4}})\varepsilon^{\frac{1}{2}}}$. Then $|L^{\varepsilon^{\frac{1}{2}}}|\leq m(\varepsilon)$.
	\end{lem}
	
	\noindent\textbf{Proof.}
	According to Lemma \ref{L1}, for any $u\in V(G)$, we deduced that
	
	\begin{eqnarray}\label{eq5}
		\sum_{v \in \overline{L_1^\eta}(u)} \sum_{w \in L_1^\eta(u) \cup L_2^\eta(u)} x_w&\leq &e(\overline{L_1^\eta}(u), L_1^\eta(u)\cup L_2^\eta(u))\nonumber\\
		&\leq &C(d_{G}(u)+|L^\eta|)\nonumber\\[2mm]
		&\leq &Cd_{G}(u)+CD(\eta,\varepsilon)\sqrt{n}
	\end{eqnarray}
	where the second inequality holds as $G$ is $\mathcal{F}$-free and $\overline{L_1^\eta}(u)\subseteq N_1(u)$. Then Lemma \ref{M5} together with (\ref{eq5}) yields
	
	\begin{equation}\label{eq6}
		\begin{aligned}
			(\beta'-1-\varepsilon^{\frac{1}{4}})nx_{u}\leq (C+x_{u})d_{G}(u)+4C\eta n.
		\end{aligned}
	\end{equation}
	
	Now, set $\eta =\varepsilon$. Then by (\ref{eq6}), for each vertex $v\in L^{\varepsilon^{\frac{1}{2}}}$, it follows that
	
	$$ d_{G}(v)\geq \frac{(\beta'-1-\varepsilon^{\frac{1}{4}})nx_{v}-4C\varepsilon n}{C+x_v} \geq \frac{(\beta'-1-2\varepsilon^{\frac{1}{4}})\varepsilon^{\frac{1}{2}}n}{2C},$$
	where the last inequality holds since $\varepsilon$ is sufficiently small and $C\geq 1$. Summing this inequality over all vertices $v\in L^{\varepsilon^{\frac{1}{2}}}$, we obtain
	
	$$|L^{\varepsilon^{\frac{1}{2}}}|\frac{(\beta'-1-2\varepsilon^{\frac{1}{4}})\varepsilon^{\frac{1}{2}}n}{2C}\leq  \sum_{v\in L^{\varepsilon^{\frac{1}{2}}}} d_{G}(v)\leq 2e(G)\leq 2Cn.$$
	Hence,
	$|L^{\varepsilon^{\frac{1}{2}}}|\leq \frac{4C^2}{(\beta'-1-2\varepsilon^{\frac{1}{4}})\varepsilon^{\frac{1}{2}}}=m(\varepsilon).$ \qed
	
	Let $u$ be an arbitrary vertex in $V(G)$. For simplicity, we write $L'=L^{\varepsilon^{\frac{1}{8}}}$, $L=L^{\varepsilon^{\frac{1}{2}}}$, $L_{i}(u)=N_{i}(u)\cap L^{\varepsilon^{\frac{1}{2}}}$, $\overline{L_{i}}(u)=N_{i}(u)\backslash L^{\varepsilon^{\frac{1}{2}}}$. Now we give a lower bound for degrees of vertices in $L'$.
	
	\begin{lem}\label{L3}
		$d_{G}(u)> (x_{u}-2\varepsilon^{\frac{1}{4}})n$ for each vertex $u\in L'$.
	\end{lem}
	
	\noindent\textbf{Proof.}
	Let $u$ be a vertex in $L'$, $S$ be a subset of $\overline{L_{1}}(u)$ such that each vertex in $S$ has at least $\beta'-1$ neighbors in $L_{1}(u)\cup L_{2}(u)$. We first prove that $|S|\leq \sqrt{n}$.  If $|L_{1}(u)\cup L_{2}(u)|\leq \beta'-2$, then $S=\emptyset$, and the inequality holds trivially. So we consider the case $|L_{1}(u)\cup L_{2}(u)|\geq \beta'-1$. Assume for contradiction that $|S|>\sqrt{n}$.
	Note that each vertex in $S$ must choose at least $\beta'-1$ neighbors from $L_{1}(u)\cup L_{2}(u)$. The number of distinct ways to choose $\beta'-1$ vertices from $L_{1}(u)\cup L_{2}(u)$ is  $\binom{|L_{1}(u)\cup L_{2}(u)|}{\beta'-1}$.
	By the pigeonhole principle, there exists a fixed subset $T\subseteq L_{1}(u)\cup L_{2}(u)$ with $|T|=\beta'-1$ such that at least $|S|/\binom{|L_{1}(u)\cup L_{2}(u)|}{\beta'-1}>\sqrt{n}/\binom{|L|}{\beta'-1} >l$ vertices in $S$ have all their $\beta'-1$ neighbors in $T$. Recall that $F_0$ is a bipartite graph in $\mathcal{F}$
	with  $|F_0|\leq l$ and $\beta'(F_0)=\beta'$. Since $u\notin L_{1}(u)\cup L_{2}(u)$ and $S\subseteq \overline{L_{1}}(u) \subseteq N_{1}(u)$, the set $\{u\}\cup T$ and the $l$ common neighbors in $S$ form a complete bipartite subgraph $K_{\beta', l}$ in $G$. Then $G$ contains $F_0$ as a subgraph , contradicting the fact that $G$ is  $\mathcal{F}$-free.
	Hence, $|S|\leq \sqrt{n}$. It follows that
	\begin{eqnarray}\label{eq7}
		\sum_{v\in \overline{L_{1}}(u)}\sum_{w\in L_{1}(u)\cup L_{2}(u)}x_{w}&\leq &e(\overline{L_{1}}(u),L_{1}(u)\cup L_{2}(u))\nonumber\\
		&=&e(\overline{L_{1}}(u)\backslash S,L_{1}(u)\cup L_{2}(u))+ e(S,L_{1}(u)\cup L_{2}(u))\nonumber\\[2mm]
		&\leq& (\beta'-2)|\overline{L_{1}}(u)\backslash S|+\sqrt{n}|L|\nonumber\\[2mm]
		&\leq& (\beta'-2)d_{G}(u)+m(\varepsilon)\sqrt{n},
	\end{eqnarray}		
		where the second to last inequality holds as $|S|\leq \sqrt{n}$ and the last inequality holds as $\overline{L_{1}}(u) \subseteq N_{G}(u)$ and $|L|\leq m(\varepsilon)$. Now, substituting $\eta=\varepsilon^{\frac{1}{2}}$ in Lemma \ref{M5} and combining with (\ref{eq7}), for large enough $n$, we derive
	\begin{eqnarray*}
		(\beta'-1-\varepsilon^{\frac{1}{4}})x_{u} n
		&\leq &d_{G}(u)x_{u}+(\beta'-2)d_{G}(u)+m(\varepsilon)\sqrt{n}+3C\varepsilon^{\frac{1}{2}} n\\
		&<&(\beta'-1)d_{G}(u)+4C\varepsilon^{\frac{1}{2}} n.
	\end{eqnarray*}
	Then $$d_{G}(u)>x_{u}n-\frac{\varepsilon^{\frac{1}{4}}}{\beta'-1}x_{u}n-\frac{4C\varepsilon^{\frac{1}{2}} }{\beta'-1}n>(x_{u}-2\varepsilon^{\frac{1}{4}})n,$$ where the last inequality holds as $\beta'\geq 2$ and $\varepsilon$ is sufficiently small. \qed
	
	\begin{lem}\label{L4}
		$(\beta'-1-2\varepsilon^{\frac{1}{4}})n\leq e(\overline{L_1}(u^{*}),\{u^{*}\}\cup L_1(u^{*}) \cup L_2(u^{*}))\leq (\beta'-1+2\varepsilon^{\frac{1}{4}})n$.
	\end{lem}
	
	\noindent\textbf{Proof.} We first prove the lower bound. By Lemma \ref{L2}, we obtain
	\begin{eqnarray}\label{eq8}
		d_{G}(u^{*})&=&d_{\overline{L_{1}}(u^{*})}(u^{*})+d_{L_{1}(u^{*})}(u^{*})\nonumber\\[2mm]
		&\leq& d_{\overline{L_{1}}(u^{*})}(u^{*})+|L_1(u^{*})|\leq d_{\overline{L_{1}}(u^{*})}(u^{*})+m(\varepsilon).
	\end{eqnarray}
	Recall that $x_{u^{*}}=1$. Setting $u=u^{*}$, $\eta=\varepsilon^{\frac{1}{2}}$ in Lemma \ref{M5}, it follows that
	
	\begin{equation}
		\begin{aligned}
			\nonumber
			(\beta'-1-\varepsilon^{\frac{1}{4}})n
			&\leq d_{G}(u^{*})+e(\overline{L_1}(u^{*}),L_1(u^{*}) \cup L_2(u^{*}))+3C\varepsilon^{\frac{1}{2}} n\\[2mm]
			&\leq e(\overline{L_1}(u^{*}),\{u^{*}\}\cup L_1(u^{*}) \cup L_2(u^{*}))+4C\varepsilon^{\frac{1}{2}} n,
		\end{aligned}
	\end{equation}
	where the last inequality holds as (\ref{eq8}) and $n$ is sufficiently large. This together with the fact that $\varepsilon$ is sufficiently small yields that
	$(\beta'-1-2\varepsilon^{\frac{1}{4}})n\leq e(\overline{L_1}(u^{*}),\{u^{*}\}\cup L_1(u^{*}) \cup L_2(u^{*}))$, as desired.
	
	Next, we prove the upper bound by contradiction. If $e(\overline{L_1}(u^{*}),\{u^{*}\}\cup L_1(u^{*}) \cup L_2(u^{*}))> (\beta'-1+2\varepsilon^{\frac{1}{4}})n$, then
	\begin{equation}\label{eq9}
		\begin{aligned}
			e(\overline{L_1}(u^{*}),L_1(u^{*}) \cup L_2(u^{*}))> (\beta'-2+2\varepsilon^{\frac{1}{4}})n.
		\end{aligned}
	\end{equation}
	
	Let $B$ be a subset of $\overline{L_1}(u^{*})$ such that each vertex in $B$ has at least $\beta'-1$ neighbors in $L_{1}(u^{*})\cup L_{2}(u^{*})$. Using the similar discussion as in the proof of Lemma \ref{L3}, we can prove that $|B|\leq \sqrt{n}$.  Combining this with Lemma \ref{L2}, we get
	\begin{eqnarray*}
		e(\overline{L_1}(u^{*}),L_1(u^{*}) \cup L_2(u^{*}))&\leq& |L|\sqrt{n}+(\beta'-2)|\overline{L_{1}}(u^{*})\backslash B|\\
		&\leq& m(\varepsilon)\sqrt{n}+(\beta'-2)n<(\beta'-2+2\varepsilon^{\frac{1}{4}})n,
	\end{eqnarray*}
	where the last inequality holds as $n$ is sufficiently large, which contradicts (\ref{eq9}).\qed
	
	\begin{lem}\label{L5}
		For all $u\in L'$,  $x_{u}\geq  1-\frac{1}{(\beta')^{3}}$ and $d_{G}(u)> \big(1-\frac{1}{(\beta')^{3}}-2\varepsilon^{\frac{1}{4}}\big)n$.
	\end{lem}
	
	\noindent\textbf{Proof.} We prove by contradiction. Suppose to the contrary that there exists a vertex $u_{0}\in L'$ such that $x_{u_{0}}< 1-\frac{1}{(\beta')^{3}}$. Recall that $L'=L^{\varepsilon^{\frac{1}{8}}}$. By Lemma \ref{L3}, $d_{G}(u_0)>(\varepsilon^{\frac{1}{8}}-2\varepsilon^{\frac{1}{4}})n$ and $d_{G}(u^{*})>(1-2\varepsilon^{\frac{1}{4}})n$. Thus,
	\begin{eqnarray}\label{eq10}
		d_{N_{1}(u^{*})}(u_0)\geq (\varepsilon^{\frac{1}{8}}-4\varepsilon^{\frac{1}{4}})n.
	\end{eqnarray}
	Then substituting $u=u^{*}$
	and $\eta=\varepsilon^{\frac{1}{2}}$ into Lemma \ref{M5}, we deduce
	\begin{eqnarray}\label{eq11}
		(\beta'-1-\varepsilon^{\frac{1}{4}})n
		&\leq & d_{G}(u^{*})x_{u^{*}}+\sum_{v \in \overline{L_1}(u^{*})} \sum_{w \in L_1(u^{*}) \cup L_2(u^{*})} x_w+3C\varepsilon^{\frac{1}{2}} n\nonumber\\
		&\leq &d_{G}(u^{*})+\sum_{v \in \overline{L_1}(u^{*})} \sum_{w \in L_1(u^{*}) \cup L_2(u^{*})\backslash\{u_0\}} x_w+d_{N_{1}(u^{*})}(u_0)x_{u_{0}}+3C\varepsilon^{\frac{1}{2}} n\nonumber\\[2mm]
		&\leq &d_{G}(u^{*})+e(\overline{L_1}(u^{*}),L_1(u^{*}) \cup L_2(u^{*})\backslash\{u_0\})+d_{N_{1}(u^{*})}(u_0)x_{u_{0}}+3C\varepsilon^{\frac{1}{2}} n\nonumber\\[2mm]
		&< &d_{G}(u^{*})+ e(\overline{L_1}(u^{*}), L_1(u^{*}) \cup L_2(u^{*}))-d_{\overline{L_1}(u^{*})}(u_0)\nonumber\\
		& & + d_{N_{1}(u^{*})}(u_0)\Big(1-\frac{1}{(\beta')^3}\Big)+3C\varepsilon^{\frac{1}{2}} n.
	\end{eqnarray}
	By  Lemma \ref{L4} and (\ref{eq8}), we obtain
	\begin{eqnarray}\label{eq12}
		& &d_{G}(u^{*})+ e(\overline{L_1}(u^{*}), L_1(u^{*}) \cup L_2(u^{*}))\nonumber\\[2mm]
		&\leq & e(\overline{L_1}(u^{*}),\{u^{*}\}\cup L_1(u^{*}) \cup L_2(u^{*}))+ m(\varepsilon)\nonumber\\[2mm]
		&\leq &(\beta'-1+2\varepsilon^{\frac{1}{4}})n+ m(\varepsilon).
	\end{eqnarray}
	Note that $d_{N_{1}(u^{*})}(u_0)-d_{\overline{L_1}(u^{*})}(u_0)=d_{L_1(u^{*})}(u_0)\leq m(\varepsilon)$. Combining  (\ref{eq10})-(\ref{eq12}), we get
	
	\begin{equation}
		\begin{aligned}
			\nonumber
			(\beta'-1-\varepsilon^{\frac{1}{4}})n
			&<  (\beta'-1+2\varepsilon^{\frac{1}{4}})n+ 2m(\varepsilon)+3C\varepsilon^{\frac{1}{2}} n-\frac{d_{N_{1}(u^{*})}(u_0)}{(\beta')^3}\\
			&<(\beta'-1+2\varepsilon^{\frac{1}{4}})n+4C\varepsilon^{\frac{1}{2}} n-\frac{(\varepsilon^{\frac{1}{8}}-4\varepsilon^{\frac{1}{4}})n}{(\beta')^3}\\
			&<\big(\beta'-1-\frac{\varepsilon^{\frac{1}{8}}}{2(\beta')^3}\big)n<(\beta'-1-\varepsilon^{\frac{1}{4}})n,
		\end{aligned}
	\end{equation}
	a contradiction. So $x_{u}\geq 1-\frac{1}{(\beta')^{3}}$ for any $u\in L'$. By Lemma \ref{L3},
	$d_{G}(u)> (1-\frac{1}{(\beta')^{3}}-2\varepsilon^{\frac{1}{4}})n$ for any $u\in L'$. \qed
	
	\vskip 2mm
	Let $\overline{L'}=V(G)\backslash L'$. Notice that $L'\subseteq L$. Thus, $|L'|\leq |L|\leq m(\varepsilon)$ by Lemma \ref{L2}.
	We will prove that $|L'|=\beta'-1$. To show it, we first give a lower bound on $e(L',\overline{L'})$.
	Since $L'\cup \overline{L'}=V(G)$, we get
	\begin{eqnarray}\label{eq13}
		\sum_{v\in L'}d_{N_{G}(u^{*})}(v)x_{v}&=&\sum_{v\in L'}d_{N_{L'}(u^{*})}(v)x_{v}+\sum_{v\in L'}d_{N_{\overline{L'}}(u^{*})}(v)x_{v}\nonumber\\
		&\leq &2e(L')+\sum_{v\in L'}d_{N_{\overline{L'}}(u^{*})}(v)x_{v}\nonumber\\
		&\leq &(m(\varepsilon))^{2}+\sum_{v\in L'}d_{N_{\overline{L'}}(u^{*})}(v)x_{v}.
	\end{eqnarray}
	On the other hand, since $e(G)\leq Cn$, we obtain
	\begin{equation}\label{eq14}
		\begin{aligned}
			\sum_{v\in \overline{L'}}d_{N_{G}(u^{*})}(v)x_{v}\leq \sum_{v\in \overline{L'}}d_{N_{G}(u^{*})}(v)\varepsilon^{\frac{1}{8}}\leq 2e(G)\varepsilon^{\frac{1}{8}}\leq 2C\varepsilon^{\frac{1}{8}}n.
		\end{aligned}
	\end{equation}
	Note that  $\lambda^{2}(G)x_{u^*}=\sum_{v\in V(G)}d_{N_{G}(u^*)}(v)x_{v}$. Summing up (\ref{eq13}) and (\ref{eq14}) gives that
	
	\begin{equation}
		\begin{aligned}
			\nonumber
			(\beta'-1-\varepsilon^{\frac{1}{4}})n &\leq \lambda^{2}(G)x_{u^{*}}=\sum_{v\in V(G)}d_{N_{G}(u^{*})}(v)x_{v}\\
			&\leq 2C\varepsilon^{\frac{1}{8}}n+(m(\varepsilon))^{2}+\sum_{v\in L'}d_{N_{\overline{L'}}(u^{*})}(v)x_{v}\\
			&<3C\varepsilon^{\frac{1}{8}}n+\sum_{v\in L'}d_{N_{\overline{L'}}(u^{*})}(v)x_{v}.
		\end{aligned}
	\end{equation}
	Then, it follows that
	
	\begin{equation}\label{eq15}
		\begin{aligned}
			e(L',\overline{L'})\geq \sum_{v\in L'}d_{N_{\overline{L'}}(u^{*})}(v)x_{v} > (\beta'-1-4C\varepsilon^{\frac{1}{8}})n.
		\end{aligned}
	\end{equation}

	\begin{lem}\label{L3.9}
		$|L'|=\beta'-1$.
	\end{lem}
	\noindent\textbf{Proof.} If $|L'|\leq \beta'-2$, then
	$$e(L',\overline{L'})\leq (\beta'-2)(n+2-\beta')<(\beta'-1-4C\varepsilon^{\frac{1}{8}})n,$$
	which contradicts (\ref{eq15}). So $|L'|\geq \beta'-1$. On the other hand, if $|L'|\geq \beta'$, then there exists a vertex set $W\subseteq L'$ with $|W|=\beta'$. By Lemma \ref{L5}, for each vertex $u\in W$,
	$d_{G}(u)> (1-\frac{1}{(\beta')^{3}}-2\varepsilon^{\frac{1}{4}})n$. Thus, $|\cap_{u\in W}d_{G}(u)|\geq (1-\frac{1}{(\beta')^{2}}-2\beta'\varepsilon^{\frac{1}{4}})n>l$, which implies that $F_{0}$ is a subgraph of $G$. This contradicts the fact that $G$ is $\mathcal{F}$-free. Hence, $|L'|=\beta'-1$. \qed
	
	\begin{lem}\label{L3.10}
		$x_{u}\geq 1-\varepsilon^{\frac{1}{10}}$ for each $u\in L'$.
	\end{lem}
	\noindent\textbf{Proof.} Suppose to the contrary that there exists a vertex $u_1\in L'$ such that $x_{u_1}< 1-\varepsilon^{\frac{1}{10}}$. By Lemmas \ref{L5} and \ref{L3.9}, we deduce that $$d_{\overline{L'}}(u_1)> (1-\frac{1}{(\beta')^{3}}-2\varepsilon^{\frac{1}{4}})n-\beta'+1>(1-\frac{2}{(\beta')^{3}})n.$$
	Then
	\begin{eqnarray}\label{eq16}
		d_{N_{\overline{L'}}(u^{*})}(u_1)x_{u_1}< d_{\overline{L'}}(u_1)(1-\varepsilon^{\frac{1}{10}})< d_{\overline{L'}}(u_1)-(1-\frac{2}{(\beta')^{3}})\varepsilon^{\frac{1}{10}}n.
	\end{eqnarray}
	On the other hand,
	\begin{eqnarray}\label{eq17}
		\sum_{v\in L'\backslash\{u_1\}}d_{N_{\overline{L'}}(u^{*})}(v)x_{v}\leq e(L'\backslash\{u_1\},\overline{L'})\leq (\beta'-1)(n-\beta'+1)-d_{\overline{L'}}(u_1).
	\end{eqnarray}
	Summing up (\ref{eq16}) and (\ref{eq17}) yields
	\begin{eqnarray*}	
		\sum_{v\in L'}d_{N_{\overline{L'}}(u^{*})}(v)x_{v}&<& (\beta'-1)(n-\beta'+1)-(1-\frac{2}{(\beta')^{3}})\varepsilon^{\frac{1}{10}}n\\
		&<&(\beta'-1-\frac{\varepsilon^{\frac{1}{10}}}{2})n<(\beta'-1-4C\varepsilon^{\frac{1}{8}})n,
	\end{eqnarray*}
	which contradicts (15). Hence,
	$x_{u}\geq 1-\varepsilon^{\frac{1}{10}}$ for each $u\in L'$, as claimed. \qed
	
	By Lemmas \ref{L3} and \ref{L3.10},  it follows that $d_{G}(u)>(1-\varepsilon^{\frac{1}{10}}-2\varepsilon^{\frac{1}{4}})n>(1-2\varepsilon^{\frac{1}{10}})n$ for each $u\in L'$. Set $R=\cap_{u\in L'}N_{\overline{L'}}(u)$. Then $|R|> (n-\beta'+1)-2(\beta'-1)\varepsilon^{\frac{1}{10}}n>(1-2\beta'\varepsilon^{\frac{1}{10}})n$.
	Take $t=|R|$. Then $K_{|L'|, |R|}=K_{\beta'-1,t}\subseteq G$.
	
	(2). By Lemma \ref{L3.10}, $x_v\geq 1-\varepsilon^{\frac{1}{10}}$ for any $v\in L'$. Recall that $L'=L^{\varepsilon^{\frac{1}{8}}}$.
	So $x_v\leq \varepsilon^{\frac{1}{8}}$ for any $v\notin L'$.
	The proof  is finished. \qed

	\begin{lem}\label{L8}
		Let $C\geq1$ be a real number,  $\mathcal{F}$ be a finite degenerate  family of graphs with $\beta'(\mathcal{F})\geq 2$,
		and $\operatorname{ex}(n, \mathcal{F})\leq Cn$. If $G$ is a graph in $\mathrm{Ex_{sp}}(n,\mathcal{F})$, then $K_{\beta'(\mathcal{F})-1,n+1-\beta'(\mathcal{F})}$ is a subgraph of $G$. Let $\mathbf{x}$ be a Perron eigenvector of $A(G)$ corresponding to $\lambda(G)$, with maximum entry equal to $1$.
		Let $L'$ be the color class of $K_{\beta'(\mathcal{F})-1, n+1-\beta'(\mathcal{F})}$ with size $\beta'(\mathcal{F})-1$. Then $x_v\geq 1-\varepsilon^{\frac{1}{10}}$ for any $v\in L'$ and  $x_v\leq \varepsilon^{\frac{1}{8}}$ for any $v\notin L'$.
	\end{lem}
	\noindent\textbf{Proof.}
	Since $K_{\beta'(\mathcal{F})-1,n+1-\beta'(\mathcal{F})}$ is $\mathcal{F}$-free, $\lambda(G)\geq  \sqrt{(\beta'(\mathcal{F})-1)(n+1-\beta'(\mathcal{F}))}$.
	Let $F_{0}\in \mathcal{F}$ be a bipartite graph with $\beta'(F_{0})=\beta'(\mathcal{F})$ and $|F_0|\leq l$.
	Define $0<\varepsilon<\big(\frac{1}{2lkC\beta'(\mathcal{F})}\big)^{120}$ be a sufficiently small constant. By Theorem \ref{T1}, $K_{\beta'(\mathcal{F})-1,t}\subseteq G$, where $t>(1-2\beta'(\mathcal{F})\varepsilon^{\frac{1}{10}})n$, and   $x_v\geq 1-\varepsilon^{\frac{1}{10}}$ for any $v\in L'$ and  $x_v\leq \varepsilon^{\frac{1}{8}}$ for any $v\notin L'$.   Set $R=\cap_{v\in L'}N_{G}(v), W=V(G)\setminus (R\cup L')$. Clearly, $|W|\leq 2\beta'(\mathcal{F})\varepsilon^{\frac{1}{10}}n$. For brevity, we write  $\beta'=\beta'(\mathcal{F})$  and $\beta=\beta(\mathcal{F})$.
	
	We claim that	
	$G$ is connected.	Suppose to the contrary that $G$ is not connected. Without loss of generality, we assume $G=G_{1}\cup G_{2}\cup\dots\cup G_p$ and $K_{\beta'-1,t}\subseteq G_{1}$. Then $G_{2}\cup\dots\cup G_p\subseteq G[W]$. Moreover, $$\lambda(G_1)\geq \lambda(K_{\beta'-1,t})> \sqrt{(\beta'-1)(1-2\beta'\varepsilon^{\frac{1}{10}})n}.$$ Recall that $|W|\leq 2\beta'\varepsilon^{\frac{1}{10}}n$ and $G[W]$ is $\mathcal{F}$-free. Then
	$$\lambda^2(G[W])\leq 2e(G[W])\leq 2C|W|\leq 4C\beta'\varepsilon^{\frac{1}{10}}n. $$ Since $\varepsilon$ is sufficiently small, we have
	$$\lambda(G_i)\leq \lambda(G[W])\leq \sqrt{ 4C\beta'\varepsilon^{\frac{1}{10}}n}<\sqrt{(\beta'-1)(1-2\beta'\varepsilon^{\frac{1}{10}})n}<\lambda(G_1)$$ for each $2\leq i\leq q$, implying  $\lambda(G)=\lambda(G_1)$. To derive a contradiction, we perform the following edge-switching operation.
	Arbitrarily choose a vertex $u_{0}\in V(G_2)$.
	Delete all edges incident to $u_{0}$ in $G$ and add all edges between $u_{0}$ and $L'$. Denote the resulting graph as $G'$. We claim $G'$ is $\mathcal{F}$-free.
	Suppose, for contradiction, that
	$G'$ contains a graph $F'\in \mathcal{F}$. Since $|R|\geq (1-2\beta'\varepsilon^{\frac{1}{10}})n$ and $\mathcal{F}$ is finite, there exists a vertex  $u_{1}\in R\setminus V(F')$. Replacing $u_0$ with $u_1$ in $F'$ yields a copy of $F'$ in $G$, which is a contradiction.

	Next we show that $W=\emptyset$.
	Suppose to the contrary that  $W\neq \emptyset$.
	Recall that $\beta'(F_{0})=\beta'$, $|F_0|\leq l$ and $F_{0}\in \mathcal{F}$.
	For any vertex  $u\in W$, we must have $d_{R}(u)\leq l-1$. Otherwise, $F_{0}\subseteq K_{\beta',l}\subseteq G$, a contradiction.
	Since $G[W]$ is $\mathcal{F}$-free,  $e(G[W])\leq C|W|$. Then there  exists a vertex $u_{0}\in W$ with $d_{W}(u_{0})\leq 2C$, which implies that $d_{W\cup R}(u_{0})\leq 2C+l-1$. Then we construct a new graph $G'$ with $V(G')=V(G)$ and $E(G')=E(G\setminus\{u_{0}\})\cup \{u_{0}u|u\in L'\}$. We claim that $G'$ is $\mathcal{F}$-free. Otherwise,  $G'$ contains a graph $F'$ in $\mathcal{F}$. Since $F'$ is finite, there exists a vertex  $u_{1}\in R\setminus V(F')$.
	Replacing $u_0$ with $u_1$ in $F'$ yields a copy of $F'$ in $G$, a contradiction. So $G'$ is $\mathcal{F}$-free. Recall that $x_{v}\leq \varepsilon^{\frac{1}{8}}$ for any  vertex $v\in V(G)\setminus L'$ and $x_{u}\geq 1-\varepsilon^{\frac{1}{10}}$ for any vertex $u\in L'$.
	Therefore,
	\begin{equation}
		\begin{aligned}
			\nonumber
			\mathbf{x}^{T}\big(\lambda(G')-\lambda(G)\big)\mathbf{x} &\geq \mathbf{x}^{T}A(G')\mathbf{x}-\mathbf{x}^{T}A(G)\mathbf{x}\\
			&\geq 2(1-\varepsilon^{\frac{1}{10}})x_{u_0}-2(2C+l-1)\varepsilon^{\frac{1}{8}}x_{u_0}\\
			&>0,
		\end{aligned}
	\end{equation}
	where the second inequality holds as $d_{L'}(u_0)\leq \beta'-2$ in $G$ and the last inequality holds as $\varepsilon$ is sufficiently small,  contradicting the maximality of
	$\lambda(G)$. So $W=\emptyset$. Thus $K_{\beta'(\mathcal{F})-1,n+1-\beta'(\mathcal{F})}\subseteq G$.
	\qed

	We present the following lemma. The core idea underlying Lemma \ref{LL2} is inspired by the work of Byrne, Desai, and Tait \cite{BDT}.

	\begin{lem}\label{LL2}
		Let $\mathcal{F}$ be a finite  family of graphs and $H=K_{\beta'(\mathcal{F})-1,n-\beta'(\mathcal{F})+1}$. For sufficiently large $n$, if $\operatorname{ex}_{H}(n,\mathcal{F})\leq e(H)+rn+O(1)$ for some $r\in [0,\frac{3}{4})$, then
		$\operatorname{ex}_{H}(n,\mathcal{F})= e(H)+rn+O(1)$, where $r\in \{0,\frac{1}{2},\frac{2}{3}\}$.
	\end{lem}

	\noindent\textbf{Proof of Lemma \ref{LL2}.}
	Let $A,B$ be the color classes of $H$ with  $|A|=\beta'(\mathcal{F})-1$ and $|B|=n+1-\beta'(\mathcal{F})$.
	Let $G$ be a graph in $\mathrm{Ex}_{H}(n,\mathcal{F})$.  We focus on the components of $G[B]$.  Assume $B'$ is the largest component of $G[B]$. We first claim $|B'|=O(1)$. Suppose to the contrary that $|B'|$ is infinite.
	Let $H'=K_{\beta'(\mathcal{F})-1,|B'|}$ and $n'=\beta'(\mathcal{F})-1+|B'|$. By the assumption of the lemma, we have
	$\operatorname{ex}_{H'}(n',\mathcal{F})< e(H')+\frac{3}{4}n'+O(1)$.
	Since $H'\subseteq G[A\cup V(B')]$,
	$$e(G[A\cup V(B')])\geq e(H')+|B'|-1=e(H')+n'-O(1)>\operatorname{ex}_{H'}(n',\mathcal{F}),$$
	contradicting the fact that  $G$ is $\mathcal{F}$-free. Next, let $B''$ be an arbitrary component in $G[B]$. If $|B''|\geq 4$, then $e(B'')\geq |B''|-1\geq  \frac{3}{4}|B''|$, which implies that the number of components of order at least 4 is $O(1)$.  Similarly, the number of components isomorphic to
	$K_3$  in $G[B]$ is $O(1)$. Define $$p=\max\{i|\  \text{the number of components in $G[B]$ with order $i$ is infinite} \}.$$ Then $p=1, 2,$ or $3$. If $p=1$, then $e(G[B])=O(1)$, so $\operatorname{ex}_{H}(n,\mathcal{F})= e(H)+O(1)$.

	If $p=2$, then $e(G[B])$ is infinite and the number of components in $G[B]$ with order at least 3 is $O(1)$. Moreover, we show that the number of isolated vertices in $G[B]$ is at most 1. Suppose to the contrary that there exist two isolated vertices $v_1$, $v_2$ in $G[B]$, then we add a new edge $e_1$ between $v_1$ and $v_2$. Denote the new graph as $G'$. We claim $G'$ is $\mathcal{F}$-free. Otherwise, if there is a graph $F_1\in \mathcal{F}$ such that $F_{1}\subseteq G'$, then $e_1 \in E(F_1)$. Combining with the fact that $F_1$ is finite and $e(G[B])$ is infinite, we always can find an edge $e_2 \in E(G[B])\setminus E(F_1)$. Thus $(F_1 \setminus{e_1})\cup e_2$ is a copy of $F_1$ in $G$, which is a contradiction. So $G'$ is $\mathcal{F}$-free. Clearly, $e(G')>e(G)$. This contradicts the fact $G\in \mathrm{Ex}_{H}(n,\mathcal{F})$. Therefore $\operatorname{ex}_{H}(n,\mathcal{F})= e(H)+\frac{1}{2}n+O(1)$.
	
	If $p=3$, then the number of $P_3$ components in $G[B]$ is infinite.  Similarly, we  can prove that the number of isolated vertices in $G[B]$ is at most 1. Next we show that the number of $P_2$ components in $G[B]$ is at most $2$. If there exist $3P_2$ components in $G[B]$, then we  delete one edge and add two new edges to get $2P_3$. Using the similar argument above, we can prove that the resulting graph is still $\mathcal{F}$-free. Moreover, the number of edges increasing, which is a contradiction. Thus,  $\operatorname{ex}_{H}(n,\mathcal{F})= e(H)+\frac{2}{3}n+O(1)$. This completes the proof. \qed
	
	Let $t=\left \lfloor \frac{n+1-\beta'(\mathcal{F})}{2} \right \rfloor$, $s=\left \lfloor \frac{n+1-\beta'(\mathcal{F})}{3} \right \rfloor$, $Q_{1}=M_{t}\cup I_{n+1-\beta'(\mathcal{F})-2t}$ and $Q_{2}=s P_{3}\cup I_{n+1-\beta'(\mathcal{F})-3s}$.
	
	\begin{lem}\label{LL14}
		Suppose $\mathcal{F}$ is a finite family of graphs. Let $H=K_{\beta'(\mathcal{F})-1,n-\beta'(\mathcal{F})+1}$ and $H_1=T\vee I_{n-\beta'(\mathcal{F})+1}$, where $T\in \mathrm{Ex}(\beta'(\mathcal{F})-1, \mathcal{H}(\mathcal{F}))$.
		Suppose $G=T\vee Q\in \mathrm{Ex}_{H_1}(n,\mathcal{F})$.
		\begin{itemize}
			\item [(1)] If $\operatorname{ex}_{H}(n,\mathcal{F})= e(H)+\frac{1}{2}n+O(1)$ and $e(Q)\geq \frac{1}{2}n-O(1)$, then $Q$ can be obtained from $Q_1$  by adding and deleting $O(1)$ edges;
			\item [(2)] If $\operatorname{ex}_{H}(n,\mathcal{F})= e(H)+\frac{2}{3}n+O(1)$ and $e(Q)\geq \frac{2}{3}n-O(1)$, then $Q$ can be obtained from $Q_2$ by adding and deleting $O(1)$ edges.
		\end{itemize}
		
	\end{lem}
	\noindent\textbf{Proof of Lemma \ref{LL14}.} Suppose $G=T\vee Q\in \mathrm{Ex}_{H_1}(n,\mathcal{F})$. If $\operatorname{ex}_{H}(n,\mathcal{F})= e(H)+\frac{1}{2}n+O(1)$, then the size of any component of $Q$ is $O(1)$ and the number of components in $Q$ with order at least 3 is $O(1)$.
	This together with the assumption $e(Q)\geq \frac{1}{2}n-O(1)$ imply that $Q$ can be obtained from $Q_1$  by adding and deleting $O(1)$ edges. If $\operatorname{ex}_{H}(n,\mathcal{F})= e(H)+\frac{2}{3}n+O(1)$ and $e(Q)\geq \frac{2}{3}n-O(1)$, using the similar discussion above
	we can show $Q$ can be obtained from $Q_2$  by adding and deleting $O(1)$ edges. \qed

	\section{Proofs of Theorems \ref{T2},  \ref{T4} and  \ref{T21}}

	\subsection{Proof of Theorem \ref{T2}}
	
	We give a lemma which will be used in the proof of Theorem \ref{T2}.

	\begin{lem}\label{LL1}
		Let $\mathcal{F}$ be a finite graph family and  $H=K_{t,n-t}$ be an $\mathcal{F}$-free graph, where $1\leq t\leq n$. For sufficiently large  $n$, the following statements are equivalent:
		\begin{itemize}
			\item[(a)] $\operatorname{ex}_{H}(n,\mathcal{F})< e(H)+\lfloor \frac{n-t}{2}\rfloor$;
			\item[(b)] there exist  integers $p,q$ such that $I_{t}\vee M_{p}$ and $I_{t}\vee S_{q+1}$ are not $\mathcal{F}$-free;
			\item[(c)] there exist a constant $c$ such that $\operatorname{ex}_{H}(n,\mathcal{F})\leq  e(H)+c$.
		\end{itemize}
	\end{lem}
	\noindent\textbf{Proof.}
	$(a)\Rightarrow (b)$.
	Let $G_{1}=I_{t}\vee \left(M_{\lfloor \frac{n-t}{2} \rfloor}\cup I_{n-t-2\lfloor \frac{n-t}{2} \rfloor}\right)$. Observe that $e(G_1)=e(H)+\lfloor \frac{n-t}{2}\rfloor>\operatorname{ex}_{H}(n,\mathcal{F})$ and $H\subseteq G_1$. Thus, $G_1$ is not $\mathcal{F}$-free, so there  exists a graph $F_{1}\in \mathcal{F}$ such that $F_{1}\subseteq G_1$.
	Since $F_{1}$ is finite, its inclusion in $ G_1$
	implies that there exists an integer $p>0$ such that $F_{1} \subseteq I_{t}\vee M_{p}$. Next,
	let $G_{2}=I_{t}\vee S_{n-t}$. Clearly, $e(G_2)=e(H)+n-t-1>\operatorname{ex}_{H}(n,\mathcal{F})$. Combining this with  $H\subseteq G_2$ and the finiteness of
	$\mathcal{F}$, we  deduce that there  exist an integer $q>0$ and a graph $F_2\subseteq \mathcal{F}$ such that $F_{2}\subseteq I_{t}\vee S_{q+1}$.

	$(b)\Rightarrow (c)$. Suppose $G$ is an arbitrary graph in $\mathrm{Ex}_{H}(n,\mathcal{F})$. Let $A, B$ be two color classes of $H$ with $|A|=t$ and $|B|=n-t$. Then $e(G[A])\leq \frac{1}{2}t^2$. By the assumption $(b)$, it follows that $\nu(G[B])\leq  p-1$ and $\Delta(G[B])\leq q-1$. According to the famous result of Chv\'atal and Hanson \cite{CH}, we have $e(G[B])\leq \nu(G[B])(\Delta (G[B])+1)\leq q(p-1)$. Therefore, $e(G)\leq e(H)+c$, where $c=\frac{1}{2}t^2+q(p-1)$.
	
	Clearly, (c) leads to (a).
	\qed

	Suppose $F_0$ is a graph in $\mathcal{F}$ with $\beta'(F_0)=\beta(\mathcal{F})$ and $|F_0|\leq l$.	
	Let $G \in \mathrm{Ex_{sp}}(n, \mathcal{F})$. By Lemma \ref{L8}, $G$ contains $K_{\beta'(\mathcal{F})-1, n+1-\beta'(\mathcal{F})}$ as a subgraph. Let $L'$ and $R$ denote the color classes of this complete bipartite graph with $|L'| = \beta'(\mathcal{F})-1$ and $|R| = n+1-\beta'(\mathcal{F})$. By Lemma \ref{LL1}, there exists a constant $c$ such that $\operatorname{ex}_H(n, \mathcal{F}) \leq e(H) + c$.
	For brevity, we write $\beta' = \beta'(\mathcal{F})$ and $\beta = \beta(\mathcal{F})$.

	\begin{lem}\label{L9}
		$G[L']\in \mathrm{Ex}(\beta'-1, \mathcal{H}(\mathcal{F})).$
	\end{lem}
	
	\noindent\textbf{Proof.} It suffices to show that $e(G[L'])=\operatorname{ex}(\beta'-1, \mathcal{H}(\mathcal{F}))$. We first prove that $e(G[L'])\geq \operatorname{ex}(\beta'-1, \mathcal{H}(\mathcal{F}))$.
	Recall that $K_{\beta'-1,n+1-\beta'}\subseteq G$ and $\operatorname{ex}_{H}(n,\mathcal{F})\leq e(H)+c$. So $e(G[R])\leq c$. Suppose for contradiction that $e(G[L'])<\operatorname{ex}(\beta'-1, \mathcal{H}(\mathcal{F}))$. We  construct a new graph $G'$ as follows:   delete all edges in $G[R]$ and $G[L']$. Embed a graph $G^{*}\in \mathrm{Ex}(\beta'-1, \mathcal{H}(\mathcal{F}))$ into the vertex set  $L'$. By Lemma \ref{L0}, $G'$ is $\mathcal{F}$-free. Applying Lemma \ref{L8}, for each $u\in L'$,
	$x_{u}\geq 1-\varepsilon^{\frac{1}{10}}$,   and for each $v\in R$,  $x_{v}\leq\varepsilon^{\frac{1}{8}}$. We obtain
\begin{equation}
		\begin{aligned}
			\nonumber
			\mathbf{x}^{T}\big(\lambda(G')-\lambda(G)\big)\mathbf{x} &\geq \mathbf{x}^{T}A(G')\mathbf{x}-\mathbf{x}^{T}A(G)\mathbf{x}\\[2mm]
			&= 2\sum_{uv\in E(G^{*})}x_{u}x_{v}-2\sum_{uv\in E(G[L'])}x_{u}x_{v}-2\sum_{uv\in E(G[R])}x_{u}x_{v}\\
			&\geq 2e(G^{*})(1-\varepsilon^{\frac{1}{10}})^{2}-2e(G[L'])-2c\varepsilon^{\frac{1}{4}} \\[2mm]
			&\geq 2-2(\beta'-1)^2\varepsilon^{\frac{1}{10}}-2c\varepsilon^{\frac{1}{4}}\\
			&>0,
		\end{aligned}
	\end{equation}
	where the second to last inequality follows from $e(G^{*})\geq e(G[L'])+1$.
	So $\lambda(G')>\lambda(G)$, which is a contradiction.
	Therefore $e(G[L'])\geq \operatorname{ex}(\beta'-1, \mathcal{H}(\mathcal{F}))$.
	
	On the other hand, suppose for contradiction that $e(G[L'])\geq \operatorname{ex}(\beta'-1, \mathcal{H}(\mathcal{F}))+1$. We  claim that $\beta'>\beta$. Assume for contradiction that $\beta'=\beta$, then $\mathcal{H}(\mathcal{F})=K_{\beta'}$ and $\operatorname{ex}(\beta'-1, \mathcal{H}(\mathcal{F}))=e(K_{\beta'-1})$. Obviously, $e(G[L'])\leq e(K_{\beta'-1})$. So $e(G[L'])\leq \operatorname{ex}(\beta'-1, \mathcal{H}(\mathcal{F}))$, contradicting the assumption $e(G[L'])\geq \operatorname{ex}(\beta'-1, \mathcal{H}(\mathcal{F}))+1$.
	Therefore, $\beta'>\beta$. Then it follows that
	$\mathcal{H}(\mathcal{F})=\mathcal{M}(\mathcal{F})$. By the definition of $\mathcal{M}(\mathcal{F})$,
	there exists a graph $\widetilde{F}\in \mathcal{M}(\mathcal{F})$ such that $\widetilde{F}\subseteq G[L']$.
	Since $\mathcal{F}$ is finite,
	there  exists an $F'\in \mathcal{F}$ such that  $F'\subseteq \widetilde{F}\vee I_{n-\beta'+1}\subseteq G$. This contradicts the fact that $G$ is $\mathcal{F}$-free. So $e(G[L'])=\operatorname{ex}(\beta'-1, \mathcal{H}(\mathcal{F}))$. Therefore $G[L']\in \operatorname{Ex}(\beta'-1, \mathcal{H}(\mathcal{F})).$\qed

	\vskip 2mm
	\noindent\textbf{Proof of Theorem \ref{T2}.}
	By Lemmas \ref{L8} and  \ref{L9}, we know $G[L']\in  \mathrm{Ex}(\beta'-1, \mathcal{H}(\mathcal{F}))$ and $G[L']\vee I_{n+1-\beta'}\subseteq G$.
	It suffices to prove $G\in \mathrm{Ex}_{G[L']\vee I_{n+1-\beta'}}(n,\mathcal{F})$. Suppose for  contradiction that $G\notin  \mathrm{Ex}_{G[L']\vee I_{n+1-\beta'}}(n,\mathcal{F})$. Then there exists a graph $G'\in \mathrm{Ex}_{G[L']\vee I_{n+1-\beta'}}(n,\mathcal{F})$ such that $G'=G[L']\vee Q$ and $e(G')\geq e(G)+1$.
	Since $\operatorname{ex}_{H}(n,\mathcal{F})\leq e(H)+c$,
	we deduce $e(G[R])+1\leq e(Q)\leq c $.

	Recall that $x_{u}\geq 1-\varepsilon^{\frac{1}{10}}$ for each $u\in L'$ and $x_{v}\leq \varepsilon^{\frac{1}{8}}$ for each $v\in R$. For a vertex $v\in R$, we have
	$$\lambda(G) x_{v}\geq \sum_{u\in L'}x_{u}\geq (\beta'-1)(1-\varepsilon^{\frac{1}{10}}),$$ which implies that $x_{v}\geq \frac{(\beta'-1)(1-\varepsilon^{\frac{1}{10}})}{\lambda(G)}$. On the other hand, since $G$ is $F_0$-free,  $d_{R}(v)\leq l-1$. Then
	\begin{equation}
		\begin{aligned}
			\nonumber
			\lambda(G) x_{v}=\sum_{u\in L'}x_{u}+\sum_{u\in N_{R}(v)}x_{u}\leq \beta'-1+(l-1)\varepsilon^{\frac{1}{8}},
		\end{aligned}
	\end{equation}
	which implies that $x_{v}\leq \frac{\beta'-1+(l-1)\varepsilon^{\frac{1}{8}}}{\lambda(G)}$. Therefore,
	\begin{eqnarray*}
		\mathbf{x}^{T}\big(\lambda(G')-\lambda(G)\big)\mathbf{x} &\geq &\mathbf{x}^{T}A(G')\mathbf{x}-\mathbf{x}^{T}A(G)\mathbf{x}\\[2mm]
		&= &2\sum_{uv\in E(Q)}x_{u}x_{v}-2\sum_{uv\in E(G[R])}x_{u}x_{v}\\
		&\geq & 2e(Q)\bigg(\frac{(\beta'-1)(1-\varepsilon^{\frac{1}{10}})}{\lambda(G)}\bigg)^{2}-2e(G[R])\bigg(\frac{\beta'-1+(l-1)\varepsilon^{\frac{1}{8}}}{\lambda(G)}\bigg)^{2}\\
		&\geq& \frac{2}{\lambda(G)^2}\bigg(\big(e(Q)-e(G[R])\big)(\beta'-1)^2-3c(\beta'-1)^2\varepsilon^{\frac{1}{10}}\bigg)\\
		&>& 0,
	\end{eqnarray*}
	where the last inequality holds as $e(Q)\geq e(G[R])+1$. Thus $\lambda(G')>\lambda(G)$, which contradicts the assumption that $G$ has the maximum spectral radius among all $\mathcal{F}$-free graphs on $n$ vertices. Hence, $G\in \mathrm{Ex}_{G[L']\vee I_{n+1-\beta'}}(n,\mathcal{F})\subseteq \mathcal{G}(\mathcal{F})$.  This completes the proof. \qed
	\subsection{Proof of Theorem \ref{T4}}

	\noindent

	Let $G\in \mathrm{Ex_{sp}}(n,\mathcal{F})$.	Since the complete bipartite graph $K_{\beta'(\mathcal{F})-1,n+1-\beta'(\mathcal{F})}$ is $\mathcal{F}$-free,  the spectral radius of $G$ satisfies $\lambda(G)\geq \sqrt{(\beta'(\mathcal{F})-1)(n+1-\beta'(\mathcal{F}))}$. By Lemma \ref{L8}, $G$ must contain $K_{\beta'(\mathcal{F})-1, n+1-\beta'(\mathcal{F})}$ as a subgraph.
	Denote the color classes of this subgraph by  $L'$ and $R$,
	where $|L'| = \beta'(\mathcal{F})-1$ and $|R|=n+1-\beta'(\mathcal{F})$. Since $G$ is $\mathcal{F}$-free, the number of edges within $L'$ satisfies $e(G[L']) \leq \operatorname{ex}(\beta'(\mathcal{F})-1, \mathcal{H}(\mathcal{F}))$.
	Suppose  $F_0 $ is a graph  in $\mathcal{F}$
	with $\beta'(F_0) = \beta'(\mathcal{F})$ and $|F_0| \leq  l$. Let $\mathbf{x}$ be the Perron vector of $A(G)$ corresponding to $\lambda(G)$, normalized such that its maximum entry $x_{u^*} = 1$. For brevity, we write  $\beta'$ instead of  $\beta'(\mathcal{F})$.

	\begin{lem}\label{L12}
		For every vertex $u \in L'$,
		\[
		1 \geq x_u \geq 1 - \frac{\beta' - 1}{\lambda(G)}.
		\]
		Furthermore, for each vertex $v \in R$, we have
		\[
		x_v = \frac{\beta' - 1}{\lambda(G)} + \Theta(n^{-1}).
		\]
		
	\end{lem}

	\noindent\textbf{Proof.} By Theorem \ref{T1}, $u^{*} \in L'$. The fact $\lambda(G)=\lambda(G) x_{u^{*}}\leq \beta'-1+\sum_{v\in R}x_{v}$ implies that $\sum_{v\in R}x_{v}\geq \lambda(G)+1-\beta'$. Thus, for each $u\in L'$, we  deduce that $$ \lambda(G) x_{u}\geq \sum_{v\in R}x_{v}\geq \lambda(G)+1-\beta'.$$ Then  $1\geq  x_{u}\geq 1- \frac{\beta'-1}{\lambda(G)}$. Moreover, $$(\beta'-1)(n+1-\beta')\leq\lambda(G)^2\leq 2e(G)\leq 2Cn.$$ So $\lambda(G)^{2}=\Theta(n)$. On the one hand, for any $v\in R$,
	$$ \lambda(G) x_{v}\geq \sum_{u\in L'}x_{u}\geq (\beta'-1)\left(1- \frac{\beta'-1}{\lambda(G)}\right),$$ which implies that $x_{v}\geq \frac{\beta'-1}{\lambda(G)}-O(n^{-1})$. Let $v^{*}$ be a vertex  in $R$ with $x_{v^{*}}=\max \{x_{v}|v\in R\}$. Since $G$ is $F_0$-free, it follows that $d_{R}(v^{*})\leq l-1$. Then
	\begin{equation}
		\begin{aligned}
			\nonumber
			\lambda(G) x_{v^{*}}=\sum_{u\in L'}x_{u}+\sum_{v\in N_{R}(v^{*})}x_{v}\leq \beta'-1+(l-1)x_{v^{*}},
		\end{aligned}
	\end{equation}
	which implies that $x_{v^{*}}\leq \frac{\beta'-1}{\lambda(G)}+O(n^{-1})$. Hence, $x_{v}=\frac{\beta'-1}{\lambda(G)}+\Theta(n^{-1})$ for any  vertex $v\in R$. This completes the proof. \qed
	
	By Lemma \ref{LL2} and the assumption of the Theorem \ref{T4},
	we have $\operatorname{ex}_{H}(n,\mathcal{F})= e(H)+r n+O(1)$, where $r=\frac{1}{2}$ or $\frac{2}{3}$. Let $c=\frac{1}{8(\beta'-1)}$,  $G_1\in \mathrm{Ex}_{H_1}(n,\mathcal{F})$. We can assume $L'\cup R$ is a vertex partition of $H_1$, and $R$ is an independent set of $H_1$.
	Then, by the assumption of Theorem \ref{T4}, we have $e(G_1)\geq \operatorname{ex}_{H}(n,\mathcal{F})-c n$, which implies that $e(G_1[R])\geq (r-c)n-O(1)$.
	
	\begin{lem}\label{L13}
		\begin{itemize}
			\item [(1)] If $r=\frac{1}{2}$, then $G[R]$ can be obtained from $Q_1$ by adding and deleting $O(1)$ edges.
			\item [(2)] If $r=\frac{2}{3}$, then $G[R]$ can be obtained from $Q_2$ by adding and deleting $O(1)$ edges.
		\end{itemize}
	\end{lem}
	\noindent\textbf{Proof.} 
	(1)	If $r=\frac{1}{2}$, then $\operatorname{ex}_{H}(n,\mathcal{F})= e(H)+\frac{1}{2}n+O(1)$, which
	implies that the size of any component of $G[R]$ is $O(1)$ and
	number of components in $G[R]$ with order at least 3 is $O(1)$.
	Furthermore,  we claim $e(G[R])$ is infinite. Otherwise, assume $e(G[R])$ is $O(1)$. Then, since $e(G_1[R])\geq (\frac{1}{2}-c)n-O(1)$, by Lemma \ref{L12},
	\begin{eqnarray*}
		\mathbf{x}^{T}\big(\lambda(G_1)-\lambda(G)\big)\mathbf{x}&\geq& \mathbf{x}^{T}A(G_1)\mathbf{x}-\mathbf{x}^{T}A(G)\mathbf{x}\\[2mm]
		&\geq& 2\operatorname{ex}(\beta'-1, \mathcal{H}(\mathcal{F}))(1-O(\lambda(G)^{-1}))^{2}-2e(G[L'])\\
		&+&2e(G_1[R])\bigg(\frac{\beta'-1}{\lambda(G)}-O(n^{-1})\bigg)^2-2e(G[R])\bigg(\frac{\beta'-1}{\lambda(G)}+O(n^{-1})\bigg)^2\\
		&\geq&2\bigg(\operatorname{ex}(\beta'-1, \mathcal{H}(\mathcal{F}))-e(G[L'])\bigg)+\frac{(1-2c)(\beta'-1)^2}{2C}-O(n^{-\frac{1}{2}}) \\
		&>&0,
	\end{eqnarray*}
	where the third equality holds as $2Cn\geq\lambda(G)^2\geq (\beta'-1)(n+1-\beta')$ and  the last inequality holds as $\operatorname{ex}(\beta'-1, \mathcal{H}(\mathcal{F}))\geq e(G[L'])$ and $c<\frac{1}{2}$. So $\lambda(G_1)>\lambda(G)$. Moreover, $G_1$ is $\mathcal{F}$-free. This contradicts the assumption that $G\in \operatorname{EX_{sp}}(n, \mathcal{F})$.
	
	Next we show the number of isolated vertices in $G[R]$ is at most 1. Suppose  to the contrary that there exist two isolated $v_1$ and $v_2$ in $G[R]$. Then we connect $v_1$ and $v_2$. We denote the
	new edge as $e_1$ and the resulting graph as $G'$. If $G'$ contains a graph $F_{1}\in \mathcal{F}$, then $e_{1}\in E(F_{1})$. Since $e(G[R])$ is infinite and $F_{1}$ is finite, there exists an edge $e_2 \in E(G[R])$ such that $e_2 \notin E(F_{1})$. Replacing $e_1$ with $e_2$ in $F_1$, we obtain a copy of $F_1$ in $G$,
	which is contradiction. Thus $G'$ is $\mathcal{F}$-free. Moreover, $\lambda(G')>\lambda(G)$, contradicting  the assumption that $G\in \mathrm{Ex_{sp}}(n, \mathcal{F})$. Therefore, the number of isolated vertices in $G[R]$ is at most $1$. Consequently, $G[R]$ can be obtained from $Q_1$ by adding and deleting $O(1)$ edges.

	(2) If $r=\frac{2}{3}$,
	we have $\operatorname{ex}_{H}(n,\mathcal{F})= e(H)+\frac{2}{3}n+O(1)$, which  implies that the size of any component of $G[R]$ is $O(1)$ and
	number of components in $G[R]$ with order at least $4$ or $K_3$ components is $O(1)$.
	Furthermore, we  show that the number of $P_3$ components in $G[R]$ is infinite. Otherwise, suppose the number of $P_3$ components in $G[R]$ is $O(1)$. Then $e(G[R])\leq \frac{1}{2}n+O(1)$. Note that $e(G_1[R])\geq (\frac{2}{3}-c)n-O(1)$. Thus,
	\begin{eqnarray*}
		\mathbf{x}^{T}\big(\lambda(G_1)-\lambda(G)\big)\mathbf{x}&\geq& \mathbf{x}^{T}A(G_1)\mathbf{x}-\mathbf{x}^{T}A(G)\mathbf{x}\\[2mm]
		&\geq &2\operatorname{ex}(\beta'-1, \mathcal{H}(\mathcal{F}))(1-O(\lambda(G)^{-1}))^{2}-2e(G[L'])\\
		& +& 2e(G_1[R]) \bigg(\frac{\beta'-1}{\lambda(G)}-O(n^{-1})\bigg)^2-2e(G[R])\bigg(\frac{\beta'-1}{\lambda(G)}+O(n^{-1})\bigg)^2\\
		&\geq &2\bigg(\operatorname{ex}(\beta'-1, \mathcal{H}(\mathcal{F}))-e(G[L'])\bigg)+\frac{(1-6c)(\beta'-1)^2}{6C}-O(n^{-\frac{1}{2}}) \\
		&>& 0,
	\end{eqnarray*}
	where the second to last inequality holds as $2Cn\geq\lambda(G)^2\geq (\beta'-1)(n+1-\beta')$ and the last inequality holds as $\operatorname{ex}(\beta'-1, \mathcal{H}(\mathcal{F}))\geq e(G[L'])$ and $c<\frac{1}{6}$. So $\lambda(G_1)>\lambda(G)$. Recall that $G_1$ is $\mathcal{F}$-free.
	This is a contradiction.
	
	Next, we demonstrate that both  the number of $P_2$ components and the number of isolated vertices in $G[R]$ are  at most $1$. Suppose there are two isolated vertices. We then  add an edge between them,  denoting the resulting graph as $G''$. Analogous to the proof in the preceding case, we can confirm that $G''$ is $\mathcal{F}$-free. Combining this  with the fact that $\lambda(G'')>\lambda(G)$, we arrive at a contradiction.  Now assume there exist two $P_2$ components $\{u_1u_2\}$ and $\{v_1v_2\}$ in $G[R]$.  Without loss of generality, we assume $x_{u_1}\geq x_{v_1}$. Define $G^{*}=G-\{v_1v_2\}+
	\{u_1 v_2\}$. By Lemma \ref{WXH}, it follows that $\lambda(G^{*})>\lambda(G)$. Moreover, since the number of $P_3$ components in $G[R]$ is infinite, one can readily verify that
	$G^{*}$ is $\mathcal{F}$-free. This contradicts the condition $G\in \mathrm{Ex_{sp}}(n, \mathcal{F})$.
	Therefore, $G[R]$ can be obtained from $Q_2$ by adding and deleting $O(1)$ edges. The proof  is complete. \qed

	\begin{lem}\label{L14}
		$G[L']\in \mathrm{Ex}(\beta'-1, \mathcal{H}(\mathcal{F})).$
	\end{lem}
	
	\noindent\textbf{Proof.} It suffices to prove $e(G[L'])=\operatorname{ex}(\beta'-1, \mathcal{H}(\mathcal{F}))$. Note that $e(G[L'])\leq \operatorname{ex}(\beta'-1, \mathcal{H}(\mathcal{F}))$ as $G$ is $\mathcal{F}$-free. Assume $e(G[L'])<\operatorname{ex}(\beta'-1, \mathcal{H}(\mathcal{F}))$. If $r=\frac{1}{2}$, then $e(G[R])\leq \frac{1}{2}n+O(1)$ and $e(G_1[R])\geq (\frac{1}{2}-c)n-O(1)$. By Lemmas \ref{L12} and  \ref{L13},  we obtain
	\begin{eqnarray*}
		\mathbf{x}^{T}\big(\lambda(G_1)-\lambda(G)\big)\mathbf{x}&\geq &\mathbf{x}^{T}A(G_1)\mathbf{x}-\mathbf{x}^{T}A(G)\mathbf{x}\\
		&\geq& 2\operatorname{ex}(\beta'-1, \mathcal{H}(\mathcal{F}))(1-O(\lambda(G)^{-1}))^{2}-2e(G[L'])\\
		&+&2e(G_1[R])\bigg(\frac{\beta'-1}{\lambda(G)}-O(n^{-1})\bigg)^2-2e(G[R])\bigg(\frac{\beta'-1}{\lambda(G)}+O(n^{-1})\bigg)^2\\
		&\geq& 2\bigg(\operatorname{ex}(\beta'-1, \mathcal{H}(\mathcal{F}))-e(G[L'])\bigg)-\frac{2c(\beta'-1)n}{n+1-\beta'}-O(\lambda(G)^{-1}) \\
		&>& \frac{3}{2}-O(n^{-\frac{1}{2}})\\
		&>&0,
	\end{eqnarray*}
	where the third inequality holds as $2Cn\geq\lambda(G)^2\geq (\beta'-1)(n+1-\beta')$ and the second to last inequality holds as $\operatorname{ex}(\beta'-1, \mathcal{H}(\mathcal{F}))\geq e(G[L'])+1$ and $c=\frac{1}{8(\beta'-1)}$.
So $\lambda(G_1)>\lambda(G)$,
which contradicts the assumption of  $G$. If $r=\frac{2}{3}$, the proof is similar to the former case, and so we omit the details here. The proof  is complete. \qed
	
	\vskip 2mm
	\noindent\textbf{Proof of Theorem \ref{T4}.}
	We proceed to prove that $G$ belongs to the family $\mathrm{Ex}_{G[L']\vee I_{n+1-\beta'}}(n,\mathcal{F})$.
	By Lemma \ref{L14}, $G[L']\in \mathrm{Ex}(\beta'-1, \mathcal{H}(\mathcal{F}))$ and  $G[L']\vee I_{n+1-\beta'}\subseteq G$. Assume for contradiction that $G\notin  \mathrm{Ex}_{G[L']\vee I_{n+1-\beta'}}(n,\mathcal{F})$. Then there exists a graph $G'\in \mathrm{Ex}_{G[L']\vee I_{n+1-\beta'}}(n,\mathcal{F})$ such that $G'=G[L']\vee Q$ and $e(G')\geq e(G)+1$. Thus, $e(Q)\geq e(G[R])+1$.
	
	If $\operatorname{ex}_{H}(n,\mathcal{F})= e(H)+\frac{1}{2}n+O(1)$, then $e(Q)\geq \frac{1}{2}n-O(1)$. By Lemma \ref{LL14}, it follows that $Q$ can be obtained from $Q_1$   by adding and deleting $O(1)$ edges.
	Combining this with Lemma \ref{L13}, we conclude that $Q$ can be obtained from  $G[R]$ by adding $p$ edges and deleting $q$ edges, where $p,q$ are two numbers satisfying $p\geq q+1$. Then, we obtain
	\begin{equation}
		\begin{aligned}
			\nonumber
			\mathbf{x}^{T}\big(\lambda(G')-\lambda(G)\big)\mathbf{x} &\geq \mathbf{x}^{T}A(G')\mathbf{x}-\mathbf{x}^{T}A(G)\mathbf{x}\\
			&\geq 2p\bigg(\frac{\beta'-1}{\lambda(G)}-O(n^{-1})\bigg)^2-2q\bigg(\frac{\beta'-1}{\lambda(G)}+O(n^{-1})\bigg)^2\\
			&\geq 2(p-q)\bigg(\frac{\beta'-1}{\lambda(G)}\bigg)^2-O(\lambda(G)^{-1}n^{-1}) \\
			&>0,
		\end{aligned}
	\end{equation}
	which contradicts the assumption that $G$ has the largest spectral radius over all $\mathcal{F}$-free graphs. Therefore, $G\in \mathrm{Ex}_{G[L']\vee I_{n+1-\beta'}}(n,\mathcal{F})\subseteq \mathcal{G}(\mathcal{F})$.
	
	If $\operatorname{ex}_{H}(n,\mathcal{F})= e(H)+\frac{2}{3}n+O(1)$, then $e(Q)\geq \frac{2}{3}n-O(1)$.
	In this case, by Lemma \ref{LL14}, it follows that $Q$ can be obtained from $Q_2$   by adding and deleting $O(1)$ edges.
	Combining this with Lemma \ref{L13}, we conclude that $Q$ can be obtained from  $G[R]$ by adding $p$ edges and deleting $q$ edges, where $p,q$ are two numbers satisfying $p\geq q+1$. Using the similar discussion above, we also can get $G\in \mathrm{Ex}_{G[L']\vee I_{n+1-\beta'}}(n,\mathcal{F})\subseteq \mathcal{G}(\mathcal{F})$. This completes the proof. \qed

	\subsection{Proof of Theorem \ref{T21}}

	\noindent\textbf{Proof of Theorem \ref{T21}.}
	Let $H=K_{\beta'(\mathcal{F})-1, n+1-\beta'(\mathcal{F})}, H_1=T'\vee I_{n+1-\beta'(\mathcal{F})}$, where $T'\in \mathrm{Ex}(\beta'(\mathcal{F})-1, \mathcal{H}(\mathcal{F}))$. By the assumption that $\{T\vee \widetilde{H}| \ T\in \mathrm{Ex}(\beta'(\mathcal{F})-1, \mathcal{H}(\mathcal{F}))\}\subseteq \mathrm{Ex}(n,\mathcal{F})$, it follows that
	$\mathcal{G}(\mathcal{F})\subseteq \mathrm{Ex}(n,\mathcal{F})$ and
	$\operatorname{ex}_{H}(n,\mathcal{F})= \operatorname{ex}_{H_1}(n,\mathcal{F})$. Evidently, $\operatorname{ex}_{H}(n,\mathcal{F})\leq \operatorname{ex}(n,\mathcal{F})\leq  e(H)+rn+O(1)$,  where $r\in [0,\frac{3}{4})$. By Lemma \ref{LL2}, we  deduce that $\operatorname{ex}_{H}(n,\mathcal{F})= e(H)+rn+O(1)$, where  $r=0, \frac{1}{2},$ or $\frac{2}{3}$.
	If $r=0$, then $\operatorname{ex}_{H}(n,\mathcal{F})< e(H)+\left \lfloor \frac{n+1-\beta'(\mathcal{F})}{2}\right \rfloor $. By Theorem \ref{T2}, we obtain that $\mathrm{Ex_{sp}}(n,\mathcal{F})\subseteq \mathcal{G}(\mathcal{F})\subseteq\mathrm{Ex}(n,\mathcal{F})$.
	If $r=\frac{1}{2}$, using Theorem \ref{T2} or Theorem \ref{T4}, we have
	$\mathrm{Ex_{sp}}(n,\mathcal{F})\subseteq \mathcal{G}(\mathcal{F})\subseteq\mathrm{Ex}(n,\mathcal{F})$. If $r=\frac{2}{3}$, Theorem \ref{T4} yields $\mathrm{Ex_{sp}}(n,\mathcal{F})\subseteq \mathcal{G}(\mathcal{F})\subseteq\mathrm{Ex}(n,\mathcal{F})$. This completes the proof.
	\qed

	\section{Proof of Theorem \ref{T3}}
	
	Let $G\in \mathrm{Ex_{sp}}(n,\mathcal{F})$. Erd\H{o}s and Gallai \cite{EG} proved that $\operatorname{ex}(n,C_{\geq k})\leq \frac{k-1}{2}(n-1)$ for any $k\geq3$. Then we have $\operatorname{ex}(n,\mathcal{F})\leq kn$.
	Since $K_{\lfloor \frac{k-1}{2} \rfloor, \ n-\lfloor \frac{k-1}{2} \rfloor}$ is $\mathcal{F}$-free, it follows that $\lambda(G)\geq \sqrt{\lfloor \frac{k-1}{2} \rfloor(n-\lfloor \frac{k-1}{2} \rfloor)}$. Let $0<\varepsilon<(\frac{1}{32k\beta'(\mathcal{F})})^{40}$ be a sufficiently small constant. By Theorem \ref{T1},
	there exists a $t\geq (1-(k+1)\varepsilon^{\frac{1}{10}})n$ such that $K_{\lfloor \frac{k-1}{2} \rfloor,\ t}$ is a subgraph of  $G$.
	We  define $L'$ as the color class
	of $K_{\lfloor \frac{k-1}{2} \rfloor,\ t}$ with size $\lfloor \frac{k-1}{2} \rfloor$ and $R=\cap_{v\in L'}N_{G}(v)$. Set $W=V(G)-R-L'$. If $e(G[R])\geq 2$, then we can find a cycle with length at least $k$ in $G$. Thus, we must have $e(G[R])\leq 1$. Moreover, when $k$ is odd, $e(G[R])=0$. Otherwise, there exists a cycle of length $k$ in $G$.
	Since $F$ is finite, we assume $|F|<l$.

	\begin{lem}
		$G$ is connected.
	\end{lem}
	\noindent\textbf{Proof.} Suppose to the contrary that $G$ is not connected.  Without loss of generality we assume that $G=G_{1}\cup G_{2}\cup\dots\cup G_p$ and $\lambda(G)=\lambda(G_1)$. Let $u_1$ be a vertex   with maximum degree in $G_1$. Then
	$$d_{G_1}(u_1)=\Delta(G_1)\geq \lambda(G)\geq \sqrt{\left\lfloor \frac{k-1}{2} \right\rfloor\left(n-\left\lfloor \frac{k-1}{2} \right\rfloor\right)}>l.$$
	Arbitrarily choose an arbitrary vertex $v_{1}\in V(G_2)$. Let $G'$ be the graph obtained from $G$ by deleting the edges incident to $v_{1}$ and adding a new edge $u_{1}v_{1}$. Clearly, $\lambda(G')>\lambda(G)$, which implies that $G'$ is not $\mathcal{F}$-free.
	Combining with the fact that $u_{1}v_{1}$ is a cut edge in $G'$ and $G$ is $C_{\geq k}$-free, we can find a copy of $F$ as a subgraph of $G'$, say $F_1$, and $v_1\in V(F_1)$.
	Since $d_{G_1}(u_1)=\Delta> l$, there exists a vertex $u'$ in $N_{G_1}(u_1)\setminus V(F_1)$.
	Then, $F_1\setminus\{v_1\}\cup\{u'\}$ is a copy of $F$ in $G$. This contradicts the fact $G$ is $\mathcal{F}$-free. Hence $G$ is connected.  \qed
	
	Let $\mathbf{x}=(x_1,\dots,x_{n})^{T}$ be a perron eigenvector of $A(G)$ corresponding to $\lambda(G)$. Let $u^{*}$ be a vertex in $V(G)$ with $x_{u^{*}}=\max\{x_{i}|\ i\in [n] \}=1$. By Theorem \ref{T1},  $x_{u}\geq 1-\varepsilon^{\frac{1}{10}}$ for each $u\in L'$ and $x_{v}\leq\varepsilon^{\frac{1}{8}}$ for each $v\in V(G)\backslash L'$.

	\begin{lem}\label{L10}
		$W=\emptyset$.
	\end{lem}
	
	\noindent\textbf{Proof.} Assume $W\neq \emptyset$. If there exists a vertex $u'\in W$ such that $d_{R}(u')\geq 2$, then we can derive a cycle with length at leas $k$ in $G$, which leads to a contradiction. Thus,  $d_{R}(u)\leq 1$ for every $u\in W$. Since $\operatorname{ex}(n,\mathcal{F})\leq kn$, we get

	\begin{equation}
		\begin{aligned}
			\nonumber
			\sum_{u\in W}d_{W\cup R}(u)\leq 2e(G[W])+e(W,R)\leq 2k|W|+ |W|=(2k+1)|W|.
		\end{aligned}
	\end{equation}
	Thus there must exist a vertex $u_0 \in W$ such that $d_{W\cup R}(u_0)\leq 2k+1$. We now perform the following operation on $u_0$: delete all edges between $u_0$ and $W\cup R$, and add all non-adjacent edges between $u_0$ and $L'$. Let $G_1$ be the resulting graph. Recall that $x_{v}\leq \varepsilon^{\frac{1}{8}}$ for any $v\in V(G)\backslash L'$  and $x_{u}\geq 1-\varepsilon^{\frac{1}{10}}$ for each $u\in L'$. Therefore,

	\begin{equation}
		\begin{aligned}
			\nonumber
			\mathbf{x}^{T}\big(\lambda(G_1)-\lambda(G)\big)\mathbf{x} &\geq \mathbf{x}^{T}A(G_1)\mathbf{x}-\mathbf{x}^{T}A(G)\mathbf{x}\\
			&\geq 2(1-\varepsilon^{\frac{1}{10}})x_{u_0}-2(2k+1)\varepsilon^{\frac{1}{8}}x_{u_0}\\
			&>0,
		\end{aligned}
	\end{equation}
	where the second inequality holds as $d_{L'}(u_0)\leq \lfloor \frac{k-1}{2} \rfloor-1$ in $G$ and the last inequality holds as $\varepsilon$ is sufficiently small.
	
	We claim that this operation does not create any copy of $F$. Otherwise,  $F$ is a subgraph of $G_1$ and  $u_0 \in V(F)$. Recall that $|F|<l$. Hence we can find a new vertex $v_0\in R\setminus V(F)$  such that $F\backslash\{u_0\}\cup\{v_0\}$ is a copy of $F$ in $G$. This is a contradiction. Consequently, we obtain a graph $G_1$ that is $F$-free and satisfies $\lambda(G_1)>\lambda(G)$.

	Let $W_1=W-\{u_0\}$. Since $G_{1}[W_1]\subseteq G[W]$,  $G_{1}[W_1]$ is $\mathcal{F}$-free. So  $e(G_{1}[W_1])\leq k|W_1|$. For each $u\in W_1$, $d_{R}(u)\leq 1$ remains valid. It follows that
	$$\sum_{u\in W_1}d_{W_1\cup R}(u)\leq  2e(G_{1}[W_1])+e(W_1,R)\leq (2k+1)|W_1|.$$
	Thus, there  exists a vertex  $u_1 \in W_1$ such that $d_{W_1\cup R}(u_1)\leq 2k+1$. Furthermore,  any vertex in $W_1$ has neighbors only in $L'$ and $W_1\cup R$. We perform the following operation on $u_1$: delete all edges between $u_1$ and $W_1\cup R$ and add  all non-adjacent edges between $u_1$ an $L'$.  Denote the resulting graph as $G_2$,  which  is $F$-free. Moreover,
	\begin{equation}
		\begin{aligned}
			\nonumber
			\mathbf{x}^{T}\big(\lambda(G_2)-\lambda(G)\big)\mathbf{x} &\geq \mathbf{x}^{T}A(G_2)\mathbf{x}-\mathbf{x}^{T}A(G)\mathbf{x}\\
			&\geq 2(1-\varepsilon^{\frac{1}{10}})x_{u_0}-2(2k+1)\varepsilon^{\frac{1}{8}}x_{u_0}+2(1-\varepsilon^{\frac{1}{10}})x_{u_1}-2(2k+1)\varepsilon^{\frac{1}{8}}x_{u_1}\\
			&>0,
		\end{aligned}
	\end{equation}
	where the second equality holds as $u_{0}, u_{1}\in W$ and the last inequality follows from $\varepsilon$ is sufficiently small. Thus $\lambda(G_2)>\lambda(G)$. Set $W_2=W_1-\{u_1\}$. Continue this process iteratively
	until $W_p=\emptyset$. At each step, $\lambda(G_i)>\lambda(G)$ holds, and the final graph
	$G_p$ is $F$-free.
	Recall that  $e(G[R])=0$ when $k$ is odd, and $e(G[R])\leq 1$ when $k$ is even. Thus $G_p=G[L']\vee  I_{n-\lfloor \frac{k-1}{2} \rfloor}$ when $k$ is odd, and
	$G_p=G[L']\vee \left (K_2 \cup I_{n-\lfloor \frac{k-1}{2} \rfloor-2}\right)$ or  $G_p=G[L']\vee  I_{n-\lfloor \frac{k-1}{2} \rfloor}$ when $k$ is even.
	Notably, all cycles in $G_p$ have length at most $k-1$.
	So  $G_p$ is $\mathcal{F}$-free and $\lambda(G_p)>\lambda(G)$, which leads to a contradiction.
	Hence $W=\emptyset$. \qed

	\begin{lem}\label{L11}
		$G[L']\in \mathrm{Ex}(\lfloor \frac{k-1}{2} \rfloor, \mathcal{H}(\mathcal{F}))$.
	\end{lem}
	
	\noindent\textbf{Proof.} To complete the proof, it suffices to show that $e(G[L'])=\operatorname{ex}(\beta'-1, \mathcal{H}(\mathcal{F}))$. We first show  $e(G[L'])\leq \operatorname{ex}(\lfloor \frac{k-1}{2} \rfloor, \mathcal{H}(\mathcal{F}))$. Assume for contradiction that $e(G[L'])\geq \operatorname{ex}(\lfloor \frac{k-1}{2} \rfloor, \mathcal{H}(\mathcal{F}))+1$.
	
	\noindent
	{\bf Case 1.} Suppose $\beta(\mathcal{F})=\lfloor \frac{k+1}{2} \rfloor=\beta'(\mathcal{F})$.
	In this case, $\mathcal{H}(\mathcal{F})=\left \{K_{\lfloor \frac{k+1}{2} \rfloor}\right\}$. Then  $e(G[L'])\leq e \left (K_{\lfloor \frac{k-1}{2} \rfloor}\right)=\operatorname{ex}(\lfloor \frac{k-1}{2} \rfloor, \mathcal{H}(\mathcal{F}))$, contradicting the assumption.
	
	\noindent
	{\bf Case 2.}  Assume $\beta(\mathcal{F})<\lfloor \frac{k+1}{2} \rfloor=\beta'(\mathcal{F})$.
	This inequality implying $\mathcal{H}(\mathcal{F})=\mathcal{M}(\mathcal{F})$. Note that $\beta(C_{\geq k})=\lfloor \frac{k+1}{2} \rfloor=\beta'(\mathcal{F})$, so
	$\mathcal{H}(\mathcal{F})=\mathcal{M}(\mathcal{F})=\mathcal{M}(F)$. By the assumption $e(G[L'])\geq \operatorname{ex}(\lfloor \frac{k-1}{2} \rfloor, \mathcal{H}(\mathcal{F}))+1$
	there must exist a graph $\widetilde{F}\in \mathcal{M}(F)$ such that $\widetilde{F}\subseteq G[L']$. Applying  Lemma \ref{L10}, since $F$ is finite, we deduce that $F\subseteq \widetilde{F}\vee I_{n-\lfloor \frac{k-1}{2} \rfloor}\subseteq G$. This contradicts the fact that $G$ is $\mathcal{F}$-free. Therefore, $e(G[L'])\leq \operatorname{ex}(\lfloor \frac{k-1}{2} \rfloor, \mathcal{H}(\mathcal{F}))$.

	Next we prove that $e(G[L'])\geq \operatorname{ex}(\lfloor \frac{k-1}{2} \rfloor, \mathcal{H}(\mathcal{F}))$ by contradiction. Assume, for the sake of contradiction, that $e(G[L'])<\operatorname{ex}(\lfloor \frac{k-1}{2} \rfloor, \mathcal{H}(\mathcal{F}))$. Let $G'=G^{*}\vee I_{n-\lfloor \frac{k-1}{2} \rfloor}$, where $G^{*}$ is  a  graph in $\mathrm{Ex}(\lfloor \frac{k-1}{2} \rfloor, \mathcal{H}(\mathcal{F}))$. By Lemma \ref{L0}, $G'$ is $\mathcal{F}$-free. Recall that $x_{u}\geq 1-\varepsilon^{\frac{1}{10}}$ for each $u\in L'$ and $x_{v}\leq\varepsilon^{\frac{1}{8}}$ for each $v\in R$. Since $e(G[R])\leq 1$, we have
	\begin{equation}
		\begin{aligned}
			\nonumber
			\mathbf{x}^{T}\big(\lambda(G')-\lambda(G)\big)\mathbf{x} &\geq \mathbf{x}^{T}A(G')\mathbf{x}-\mathbf{x}^{T}A(G)\mathbf{x}\\
			&= 2\sum_{uv\in E(G^{*})}x_{u}x_{v}-2\sum_{uv\in E(G[L'])}x_{u}x_{v}-2\sum_{uv\in E(G[R])}x_{u}x_{v}\\
			&\geq 2e(G^{*})(1-\varepsilon^{\frac{1}{10}})^{2}-2e(G[L'])-2\varepsilon^{\frac{1}{4}} \\
			&\geq 2-2\left (\left \lfloor \frac{k-1}{2} \right\rfloor\right)^2\varepsilon^{\frac{1}{10}}-2\varepsilon^{\frac{1}{4}}\\
			&>0,
		\end{aligned}
	\end{equation}
	where the second to last inequality follows from $e(G^{*})\geq e(G[L'])+1$, which contradicts that $G\in \mathrm{Ex_{sp}}(n,\mathcal{F})$. Therefore, $G[L']\in \mathrm{Ex}(\lfloor \frac{k-1}{2} \rfloor, \mathcal{H}(\mathcal{F}))$. \qed
	
	\vskip 2mm
	\noindent\textbf{Proof of Theorem \ref{T3}.}
	Recall that  $e(G[R])=0$ when $k$ is odd, and $e(G[R])\leq 1$ when $k$ is even.
	If $k$ is odd,  by Lemmas \ref{L10} and  \ref{L11},
	$$G=G[L']\vee I_{n- \frac{k-1}{2} }\in  \left\{T\vee I_{n- \frac{k-1}{2}}\bigg| \ T\in\mathrm{Ex}\left( \frac{k-1}{2}, \mathcal{H}(\mathcal{F})\right)\right\}. $$
	If $k$ is even, by Lemmas \ref{L10} and  \ref{L11},
	$$G=G[L']\vee I_{n- \frac{k}{2}+1}\in \left\{T\vee I_{n- \frac{k}{2}+1}\bigg|\ T\in\mathrm{Ex}\left( \frac{k}{2}-1, \mathcal{H}(\mathcal{F})\right)\right\} $$ or $$G=G[L']\vee \left(K_2 \cup  I_{n- \frac{k}{2}-1}\right)\in \left \{T\vee \left(K_2 \cup  I_{n- \frac{k}{2}-1}\right )\bigg|\ T\in\mathrm{Ex}\left( \frac{k}{2}-1, \mathcal{H}(\mathcal{F})\right)\right\}.$$ This completes the proof. \qed


\begin{thebibliography}{99}
		\bibitem{AF}
		N. Alon, P. Frankl, Tur\'an graphs with bounded matching number, J. Comb. Theory, Ser. B 165 (2024) 223-229.
		
		\bibitem{B}
		J. Byrne, A sharp spectral extremal result for general non-bipartite graphs, arXiv preprint arXiv: 2411.18637, (2024).
		
		\bibitem{BDT}
		J. Byrne, D.N. Desai, M. Tait, A general theorem in spectral extremal graph theory, arXiv preprint arXiv: 2401.07266, (2024).
		
		\bibitem{CLZ2019}		
		M. Chen, A. Liu,  X. Zhang, Spectral Extremal Results with Forbidding
		Linear Forests, Graphs and Comb. 35 (2019) 335-351.
		
		\bibitem{CLZ2021}
		M. Chen, A. Liu,  X. Zhang, On the spectral radius of graphs without
		a star forest, Discrete Math. 344(4) (2021) 112269.
		
		\bibitem{CDT}
		S. Cioab\u{a}, D.N. Desai, M. Tait, A spectral Erd\H{o}s-S\'os theorem, SIAM J. Discrete Math. 37(3) (2023) 2228-2239.
		
		\bibitem{PC2011}
		P\'eter Csikv\'ari, Applications of the Kelmans transformation: extremality of the threshold graphs, Electron. J. Comb. 18 (2011) \#P182.
		
		
		\bibitem{CDT0}
		S. Cioab\u{a}, D.N. Desai, M. Tait, The spectral radius of graphs with no odd wheels, Eur. J. Comb. 99 (2022) 103420.
		
		\bibitem{CDT}
		S. Cioab\u{a}, D.N. Desai, M. Tait, A spectral Erd\H{o}s-S\'os theorem, SIAM J. Discrete Math. 37(3) (2023) 2228-2239.
		
		\bibitem{CDT1}
		S. Cioab\u{a}, D.N. Desai, M. Tait, The spectral even cycle problem, Comb. Theory 4(1) (2024) 10.
		
		\bibitem{CH}
		V. Chv\'atal, D. Hanson, Degrees and matchings, J. Comb. Theory, Ser. B 20 (1976) 128-138.
		
		\bibitem{DHP}
		C. Dou, F. Hu, X. Peng, Tur\'an numbers of cycles plus a general graph, arXiv preprint arXiv:2411.17322, (2024).
		
		\bibitem{DNP}
		C. Dou, B. Ning, X. Peng, The number of edges in graphs with bounded clique number and circumference, arXiv preprint arXiv:2410.06449, (2024).
		
		\bibitem{EG}
		P. Erd\H{o}s, T. Gallai, On maximal paths and circuits of graphs, Acta Math. Acad. Sci. Hung. 10 (1959) 337-356.
		
		\bibitem{ES}
		P. Erd\H{o}s, M. Simonovits, A limit theorem in graph theory, Studia Sci. Math. Hungar. 1 (1966) 51-57.
		
		\bibitem{FS}
		Z. F\"uredi, M. Simonovits, The history of degenerate (bipartite) extremal graph problems, Erd\H{o}s centennial, 169-264, Bolyai Soc. Math. Stud., 25, J\'anos Bolyai Math. Soc., Budapest, 2013.
		
		\bibitem{FTZ}
		L. Fang, M. Tait, M. Zhai, Decomposition family and spectral extremal problems on non-bipartite graphs, Discrete Math. 348 (2025) 114527.
		
		\bibitem{FLSZ2024}		
		L. Fang, H. Lin, J. Shu,  Z. Zhang, Spectral extremal results on trees,
		Electron. J. Comb. 31(2) (2024), \#P2.34.
		
		\bibitem{FZC}
		X. Fang, X. Zhu, Y. Chen, Generalized Tur\'an problem for a path and a clique, Eur. J. Comb. 127 (2025) 104137
		
		\bibitem{FYZ2007}		
		L. Feng, G. Yu, and X. Zhang,  Spectral radius of graphs with given matching number, Linear Algebra Appl. 422(1) (2007) 133-138.
		
		\bibitem{GH2019}
		J. Gao and X. Hou, The spectral radius of graphs without long cycles, Linear Algebra Appl. 566 (2019) 17-33.
		
		\bibitem{G}
		D. Gerbner, On Tur\'an problems with bounded matching number, J. Graph Theory 106 (2024) 23-29.
		
		
		\bibitem{JYZ}
		S. Jiang, Y. Zhai, X. Yuan, Some stability results for spectral extremal problems of graphs with bounded matching number, Linear Algebra Appl. 708 (2025) 513-524.
		
		\bibitem{KX}
		G. Katona, C. Xiao, Extremal graphs without long paths and large cliques, Eur. J. Comb. 119 (2024) 103807.
		
		\bibitem{LK}
		Y. Liu, L. Kang, Extremal graphs without long paths and a given graph, Discrete Math. 347 (2024) 113988.
		
		\bibitem{N2010}
		V. Nikiforov,  The spectral radius of graphs without paths and cycles of specified length, Linear Algebra Appl.  432(9) (2010) 2243-2256.
		
		\bibitem{N2}
		V. Nikiforov, Merging the $A$- and $Q$-spectral theories, Appl. Anal. Discrete Math. 11(1) (2017) 81-107.
		
		\bibitem{N}
		V. Nikiforov, Some new results in extremal graph theory. Surveys in combinatorics
		2011, 141-181, London Math. Soc. Lecture Note Ser., 392, Cambridge Univ. Press,
		Cambridge, 2011.
		
		\bibitem{N1}
		V. Nikiforov, The spectral radius of graphs without paths and cycles of specified length, Linear Algebra Appl. 432 (2010) 2243-2256.
		
		\bibitem{S}
		M. Simonovits, A method for solving extremal problems in graph theory, stability problems, in: Theory of Graphs, Proc. Colloq., Tihany, 1966, Academic Press, 1968, pp. 279-319.
		
		\bibitem{WFL}
		T. Wang, L. Feng, L. Lu, Spectral extremal problems for graphs with bounded clique number, Linear Algebra Appl. 710 (2025) 273-295.
		
		\bibitem{WHM}
		H. Wang, X. Hou, Y. Ma, Spectral extrema of graphs with bounded clique number and matching number, Linear Algebra Appl. 669 (2023) 125-135.
		
		\bibitem{WKX}
		J. Wang, L. Kang, Y. Xue, On a conjecture of spectral extremal problems, J. Comb. Theory, Ser. B 159 (2023) 20-41.
		
		\bibitem{WXH}
		B. Wu, E. Xiao, Y. Hong, The spectral radius of trees on $k$ pendant vertices, Linear Algebra Appl. 395 (2005) 343-349.
		
		\bibitem{XK}
		Y. Xue, L. Kang, On generalized Tur\'an problems with bounded matching number, arXiv preprint arXiv:2410.12338, (2024).
		
		\bibitem{Y}
		L. Yuan, X. Zhang, Tur\'an numbers for disjoint paths, J. Graph Theory 98 (2021) 499-524.
		
		
		\bibitem{ZC}
		X. Zhu, Y. Chen, Extremal problems for a matching and any other graph, J. Graph Theory 109 (2025) 19-24.
		
		\bibitem{ZL}
		X. Zhao, M. Lu, Generalized Tur\'an problems for a matching and long cycles, arXiv preprint arXiv:2412.18853, (2024).
		
		\bibitem{ZY}
		Y. Zhai, X. Yuan, Spectral extrema of $\{K_{k+1}, \mathcal{L}_{s}\}$-free graphs, Linear Algebra Appl. 682 (2024) 309-322.
		
		\bibitem{ZYY}
		Y. Zhai, X. Yuan, L. You, Spectral extrema of graphs: Forbidden star-path forests, Discrete Math. 348 (2025) 114351.
		
	\end{thebibliography}
\end{document}